\documentclass[11pt]{article}
\usepackage{}
\usepackage{mathrsfs}
\usepackage{amsmath}
\usepackage{amsfonts}
\usepackage{amssymb}
\usepackage{epsfig}
\usepackage{dsfont}
\usepackage{multirow}
 \usepackage{graphicx}
\usepackage{color}

\usepackage[colorlinks=true]{hyperref}
\hypersetup{urlcolor=blue, citecolor=blue,linkcolor=blue}

\renewcommand{\Box}{\framebox{\rule{0.3em}{0.0em}}}

\def\prob {{\rm Prob}}

\newtheorem{thm}{Theorem}[section]
\newtheorem{lema}{Lemma}[section]
\newtheorem{prop}{Proposition}[section]
\newtheorem{ex}{Example}[section]

\newtheorem{defi}{Definition}[section]
\newtheorem{alg}{Algorithm}[section]
\newtheorem{assu}{Assumption}[section]

\newcommand{\setd}{{ d \kern -.15em l}}
\newcommand{\hatsetd}{ d \hat{\kern -.15em l }}
\newcommand{\dd}{\mathsf {d\kern -0.07em l}} %	nested

\newcommand{\bgeqn}{\begin{eqnarray}}
\newcommand{\edeqn}{\end{eqnarray}}
\newcommand{\bgeq}{\begin{eqnarray*}}
\newcommand{\edeq}{\end{eqnarray*}}
\newcommand{\bec}{\begin{center}}
\newcommand{\enc}{\end{center}}
\newcommand{\R}{{\rm I\!R}}

\newcommand{\inmat}[1]{\mbox{\rm {#1}}}

\newcommand{\B}{{\cal B}}

\newcommand{\be}{\begin{equation}}
\newcommand{\ee}{\end{equation}}

\def\w{\omega}
\def\bbe{{\Bbb{E}}} %expectation

\renewcommand{\Box}{\hfill \rule{2.3mm}{2.3mm}}

\begin{document}
%\textcolor{Red}{The content in red color is not resolved.}

%
%\textbf{This document is not suitable for circulation}
%rough draft  not suitable for circulation}

\begin{center}

{\LARGE %(Stability analysis)
Continuous Behavioural Function Equilibria
and Approximation Schemes in Bayesian Games with Non-Finite Type and Action Spaces
}
% %Asymptotic Convergence Analysis for Distributional Robust Optimization and Equilibrium Problems\footnote{The research is supported by EPSRC %grant EP/J014427/1.}
%% Robust games and equilibrium}
%% On the e
%%Existence and Numerical Solution of
%%equilibria to
%Bayesian
%%games
%%Equilibria
%Games with Non-Finite Type and Action Spaces:
%Existence of Continuous Behavioural Function Equilibria
%and
%%Numerical
%Approximation Schemes}
% %Stability Analysis for Mathematical Programs with Robust Chance Constraint}
% %}

\vspace{0.6cm}

Shaoyan Guo
\\
%Institute of Operations Research and Control Theory,
 School of Mathematical Sciences, Dalian University
of Technology,
 Dalian 116024, China \\(syguo@dlut.edu.cn)

Huifu Xu
\\
School of Mathematical Sciences,
   University of Southampton,
 SO17  1BJ, Southampton, UK\\
(H.Xu@soton.ac.uk)

Liwei Zhang
\\
%Institute of Operations Research and Control Theory,
 School of Mathematical Sciences, Dalian University
of Technology,
 Dalian 116024, China \\(lwzhang@dlut.edu.cn)

\vspace{0.6cm}
\today

\end{center}

%\maketitle

\noindent
\textbf{Abstract.}
Meirowitz \cite{Mei03} showed existence of continuous behavioural function equilibria
for Bayesian games with non-finite type and action spaces. A key condition for the proof of the existence result is
equi-continuity of behavioural functions which, according to  Meirowitz \cite[page 215]{Mei03},
is likely to fail or difficult to verify. In this paper, we advance the research by presenting some
verifiable conditions for the required equi-continuity,
namely
some growth conditions
%such as
%the second order growth conditions
of
the expected utility functions of each player at equilibria.
In the case when the growth is of second order, we demonstrate that
%The second order growth
the condition is guaranteed by strong concavity of the utility function.
Moreover, by using recent research on polynomial decision rules and optimal discretization approaches in stochastic and robust optimization, we propose some approximation schemes for the Bayesian equilibrium problem: first, by restricting the behavioural functions
to polynomial functions of certain order over the space of types, we demonstrate that solving a Bayesian polynomial behavioural function
equilibrium is down to solving a finite dimensional stochastic equilibrium problem; second, we apply
the optimal quantization
%approach
method due to Pflug and Pichler \cite{PfP11} to develop an effective
discretization scheme for solving the latter.
Error bounds are derived
for the respective approximation schemes under moderate conditions and both academic examples and numerical
results are presented to explain the Bayesian equilibrium problem and their approximation schemes.


\section{Introduction}

Over the past few years, there has been an increasing attention to Nash games with private information.
 A common assumption in such games is that the prior distribution of the types of all players is known in public,
 each player has complete information of its own type which determines its utility function but is unaware of its rival's type.
Based on the prior information, each player chooses a response function which is also known as behavioural function
defined over its type space under Nash conjecture
and an equilibrium arising from this kind of game is known as Bayesian Nash equilibrium, see
Hansanyi \cite{Har67} for a comprehensive original discussion of the Bayesian games where private information might also include
other aspects of a player's payoff function.

Meirowitz \cite{Mei03} considered a Bayesian game
 where each player chooses its behavioural function based
 on  maximization of its expected utility with
  the expectation being taken w.r.t. its rival's distribution
  of types conditional on the selection of its own type.
      Under some conditions, he established existence of equilibria for the
      Bayesian game using Schauder's fixed point theorem.
      One of the main conditions that Meirowitz used
      for the existence result
      is equi-continuity of the behavioural functions
      which is
      elicited to ensure that the space of behavioural functions
      is closed and the operator mapping the set of
 behaviour functions to itself
is compact.
%relatively compact set-valued.
Meirowitz commented that the equi-continuity condition is likely to fail or difficult to verify in practical applications.
%{\color{blue}
Athey \cite{At01} considered a class of Bayesian games  where the types are drawn from an atomless joint probability distribution and
each player's utility function has so-called single crossing property which means  whenever each opponent uses a nondecreasing strategy in the sense that higher types choose higher actions, a player's best response strategy is also nondecreasing.
Under these circumstances, she demonstrated existence of equilibria in every finite-action game
with each player's behavioural function being nondecreasing and step-like. Moreover, when the space of actions is continuous,
she showed existence of a sequence of nonincreasing step-like (behavioural function) equilibria to finite action games
that converges to an equilibrium with the continuum-action which means that an equilibrium in continuous action spaces can be approximated by a sequence of  nondecreasing step-like behavioural function equilibria in finite action spaces.

%%%%%%%%%%%%%%%%%%%

Ui \cite{Ui16}
provided a sufficient condition for the existence and uniqueness of a Bayesian Nash equilibrium
by regarding it as a solution of a variational inequality where the payoff gradient of the game is defined as a vector
whose component is a partial derivative of each player’s payoff function with respect to the player’s own
action. He demonstrated that when the Jacobian matrix of the payoff gradient is negative definite for each type,
a Bayesian Nash equilibrium exists using some theories in variational inequality rather than Schauder's fixed-point theorem.
Note that the Bayesian Nash equilibrium considered by Ui \cite{Ui16}
is slightly different from  Meirowitz's where a player's behavioural function is optimal almost surely
for its type. This means the behavioural function is not necessarily optimal at a subset of its type set with Lebesgue measure
 zero. In some references, this kind of equilibrium is called pure strategy Nash equilibrium (PSNE), see \cite{EC14, Fey08}.
Of course, the behavioural functions at such an equilibrium are not necessarily continuous.

A particular interesting application area of the Bayesian equilibrium model is
Tullock's rent-seeking
contest \cite{Tull80, Tul67}.
%contest.
 A rent-seeking contest is a situation where players
spend costly efforts to gain a reward. Many conflict situations can be described by  rent-seeking contests including political
campaigns, patent races, war fighting, lobbying efforts, labor market competition, legal battles and professional sports,
see Fey \cite{Fey08} and references therein. Fey showed existence of symmetric Bayesian equilibrium in the case when there are two players in the contest.
% and each player's effort is restricted to
Ewerhart \cite{EC14} advanced the research by showing existence of a unique PSNE where the contest success function is of logit form with concave impact functions and player's private information may relate to either costs or valuations.

Aghassi and Bertsimas \cite{AgB06} discussed a broad class of robust games with finite number of players,
each  player plays a mixed strategy over a finite set of pure strategies and the optimal response is based on the worst payoff matrix. In particular, they  investigated robust games with private information where each player's behavioural function is based on the worst type and worst payoff matrix. Under some conditions, they established existence of robust equilibria using a fixed point theorem due to Bohnenblust and Karlin \cite{BoK50}. A key element in the existence theorem is compactness: the set of behavioural functions must be compact and
the mapping which takes each behavioural function to a subset of the behavioural functions is compact and convex set-valued. By Arzela-Ascoli's theorem (see \cite{Mei03}), the latter compactness
is fulfilled if and only if behavioural functions in the image space are bounded and equi-continuous.

In this paper, we
%take on
extend the research in two directions.
%by extending the existing results.
%ence result.
First, we derive verifiable sufficient conditions for equi-Lipschitz continuity of the behavioural functions, a
key condition used by Meirowitz \cite{Mei03} for showing the existence of an equilibrium.
% in his Bayesian game model.
This might help to make his model and the equilibrium results more applicable.
Second, we apply the well-known decision rules for calculating an approximate behavioural function equilibrium.
The fundamental idea is to restrict the behavioural function of each player to polynomial functions
of certain order.
In doing so, we will be able to effectively converting the Bayesian game
into a finite dimensional
stochastic game model which can be solved by existing
stochastic approximation methods
such as sample average approximation method and optimal quantization method.
The approach is known as polynomial decision rules in the literature of stochastic optimization and robust optimization, see
%The idea comes from
for instances Bampou and Kuhn \cite{BK12}
%who successfully applied
for the polynomial decision rules applied to
%approximate
 continuous linear programs and Kuhn, Wiesemann and Georghiou \cite{KWG11} for linear decision rules applied to
 distributionally robust formulation of two stage stochastic programs.
%Over the past few years, linear and polynomial decision rules have been shown a very effective approach to
%approximate infinite dimensional solutions in stochastic and robust
%optimization problems, see
%\cite{BK12,KWG11} and references therein.

As far as we are concerned, the main contributions of this paper can be summarized as follows.

\begin{itemize}

\item We revisit
%the Bayesian game model
the existence results established
%considered
by Meirowitz \cite[Proposition 1]{Mei03} for the Bayesian game
by replacing the explicit assumption of equi-conituity
of the behavioural functions with
some
%the second order
growth conditions of the expected utility functions of each player at equilibria
(Theorem \ref{t-NE-existence}). The new existence result is derived by using
a general stability result in
% pragmatic
parametric programming (Lemma \ref{l-stability}).
In the case when the growth is of second order,
a sufficient condition
%s for the second order growth condition
is given (Proposition \ref{p:second-growth}).
Moreover,
%in the case
{\color{black}
when 
%the gradient of 
the utility functions of all players are directionally differentiable and satisfy 
certain monotonicity conditions,
% are strictly monotone 
with respect to their actions, 
we demonstrate uniqueness of the Bayesian equilibrium (Theorem \ref{t-unique-BNE}).
}

\item We propose to use polynomial decision rules to derive an approximation of the behavioural functions and hence
the Bayesian behavioural function equilibria. This is possible when we concentrate on the continuous Bayesian equilibrium model (Theorem \ref{p:BNE-inequality}).  Under the approximation framework, we demonstrate existence of polynomial behavioural function equilibria (Theorem \ref{t-exist-poly-BNE}) and show that solving a Bayesian polynomial behavioural function
equilibrium is down to finding a finite dimensional stochastic equilibrium problem.
Convergence of polynomial Bayesian equilibrium to the true Bayesian equilibrium is established to justify the polynomial decision rules.
 Moreover, we apply the optimal quantization approach due to Pflug and Pichler \cite{PfP11} to develop an effective
discretization scheme for solving the approximate Baygesian equilibrium model.
Error bounds are derived
for the approximation schemes under moderate conditions and both academic examples and numerical
results are presented to explain the Bayesian equilibrium problem and their approximate schemes (Theorem \ref{t-V^MtoV}).

\item We
apply the proposed theory of existence and uniqueness of behavioural function
equilibrium
and the approximation schemes
to rent-seeking contests.
Specifically,
for general symmetric
multi-player games,
we show that our conditions of existence and uniqueness can be easily satisfied
%the rent-seeking contest game where
when each player's effort is lower bounded by a positive number.
In other words, we can show existence and uniqueness of a continuous
behavioural function
equilibrium rather than a PSNE.
Moreover, by driving the lower bound to zero,
 we show that the sequence of the behavioural function
equilibria
  %a slightly stronger conclusion
 %that
has at least
% there exists
a cluster point which
% by showing $f^*$
is a continuous behavioural function equilibrium rather than
a PSNE
 of the unconstrained contest where the player's effort does not have a positive lower bound,
slightly strengthening Ewerhart's earlier result \cite[Theorem 3.4]{EC14}, see Proposition \ref{p:rent-geq0}.

\end{itemize}

%{\color{red}
%1.
%%Review on Athey's paper;
%literature review on Bayesian game;\\
%2.Review on the variational inequality in an infinite-dimensional space and how to solve it;\\
%%3.Review on polyhedral decision rule, refer to Daniel Kuhn's paper: primal and dual...2011;
%3. Benefit of polynominal approximation.
%}

The rest of the paper is organized as follows.  In Section 2, we present a detailed
explanation of the Bayesian Nash equilibrium  model, its equivalent
formulations and key difference between behavioural function
equilibrium  and so-called pure strategy Nash equilibrium.
In Section 3, we investigate existence and uniqueness
of Bayesian Nash equilibrium based on new conditions which sufficiently
ensure behavioural function of each player to be equi-Lipschitz continuous.
 In Section 4, we
 %apply
 discuss approximation schemes for the Bayesian Nash equilibrium model, we start with
   polynomial decision rules and then followed by
optimal quantization schemes, convergence results are derived to justify the approximations.
 Finally, in Section 5, we examine the established theory and approximation schemes
 by applying them to rent-seeking contests and present preliminary numerical test results.

\section{The model}

We consider a Bayesian game with $n$ players.
%{\color{blue}
Each player possesses
a preference utility function
denoted by $u_i(a_i, a_{-i},\theta_i,\theta_{-i})$ for $i=1,\cdots,n$
which  depends on
%player $i$'s
the player's
action $a_i$, its rival's actions $a_{-i}$,
the player's type $\theta_i$ and
the rival's type $\theta_{-i}$.
We assume that a type $\theta_i$ takes values from set $\Theta_i$
%\subset \R^{d_i}$
and an action $a_i$ takes values from action space $\mathcal{A}_i$
%\subset \R^{z_i}$,
where $\Theta_i$ and $\mathcal{A}_i$ are non-empty, compact and convex subsets
of $\R^{d_i}$ and $\R^{z_i}$ respectively.
%}
Following the terminology of Meirowitz \cite{Mei03}, a profile of types is a vector $\theta=(\theta_1,\cdots,\theta_n)\in \Theta:=\Theta_1\times \cdots\times \Theta_n$ and a profile of actions is a vector $a=(a_1,\cdots,a_n)\in \mathcal{A} :=\mathcal{A}_1\times \cdots\times \mathcal{A}_n$.
 Using the standard notation, we denote by $a_{-i}$ and $\theta_{-i}$
 respectively the vector of actions and the
  types  of all players except $i$.
Conditional on its type $\theta_i$,
 player $i$'s posterior belief about $\theta_{-i}$ is represented by a
 conditional probability distribution $\eta_i(\cdot|\theta_i)$,
 %{\color{blue}
 which describes the probability of player $i$'s rivals taking
 a particular type $\theta_{-i}$.
 %}
%Player's

%{\color{blue}
Information on players' types is private which means each player
only knows its own type but not other's.
However, it is assumed that
the probability distribution of $\theta$, denoted by $\eta(\theta)$, is
public information. This information describes the probability of all players taking
a particular $\theta$ which may be retrieved from empirical data.
%}
%{\color{blue}
Throughout the paper, we will use  $\theta$ to denote
 a deterministic element of $\R^{d_1+\cdots+d_n}$
 or a random vector $\theta(\w)$ mapping from probability space $(\Omega,\mathscr{B},\eta)$
to $\R^{d_1+\cdots+d_n}$
 depending on the context.
%}

For $i=1,\cdots,n$, we denote by $\mathcal{F}_i$  the set of functions $f_i:\Theta_i\to \mathcal{A}_i$ with the infinity norm,
that is
$$\|f_i\|_{\infty}=\max_{\theta_i \in \Theta_i} |f_i(\theta_i)|,
$$
and $\mathcal{C}_i$  the set of continuous functions $f_i:\Theta_i\to \mathcal{A}_i$.
% such that
%$$
%\int_{\theta_i\in \Theta_i} \|f_i(\theta_i)\|^pd\theta_i<\infty
%$$
%where $p\in (1,\infty)$ is a constant.
%Obviously $F_{i}$ forms an $\mathscr{L}^p$ space,
%{\color{red}
Equipped with the infinity norm,
%Obviously
 $\mathcal{C}_{i}$ forms a closed, bounded and convex Banach space.
%}
%We may write $F_i = \mathscr{L}^p(\Theta_i,A_i)$.
For the simplicity of notation, let
\bgeqn
\mathcal{F} :=\mathcal{F}_1\times\cdots \times \mathcal{F}_n, \quad \mathcal{C} :=\mathcal{C}_1\times\cdots \times \mathcal{C}_n,
\label{eq:F-C}
\edeqn
 %$\mathscr{L}^p(\Theta)=\mathscr{L}^p(\Theta_1,\R^{z_1})\times\cdots \times\mathscr{L}^p(\Theta_n,\R^{z_n})$
 and $N :=\{1,\ldots,n\}$.

%If for each $i =1,\ldots, n$, $A_i \subset \R^{z_i}$ is nonempty compact, then $F_i$ is a convex,
% bounded and closed set in $\mathscr{L}^p(\Theta_i,\R^{z_i})$. It follows from Corollary 3.22 of \cite{Brezis}
%that $F$ is weakly compact in $\mathscr{L}^p(\Theta)$.

\begin{defi}[Bayesian behavioural function equilibria]\rm
{\color{black}A  {\em behavioural function equilibrium}}
is an
$n$-tuple  $f=(f_1,\ldots,f_n)$
mapping from $\Theta_1\times\cdots\times \Theta_n$ to
$\mathcal{A}_1\times \cdots \times \mathcal{A}_n$
such that for every $i\in N$,
\bgeqn
\label{eq:BNE}
\inmat{(BNE)}  \quad f_i(\theta_i) \in\arg \max_{a_i\in \mathcal{A}_i}
\int_{\theta_{-i}\in \Theta_{-i}} u_i(a_i, f_{-i}(\theta_{-i}),\theta_i,\theta_{-i})d\eta_i
(\theta_{-i}|\theta_i), \forall \theta_i\in \Theta_i,
\edeqn
where
$
\Theta_{-i} :=
\Theta_1\times \cdots \times \Theta_{i-1}\times\Theta_{i+1}\times \cdots \times \Theta_n$,
$\eta_i(\theta_{-i}|\theta_i)$ is the conditional probability distribution
of $\theta_{-i}$, that is,
$\eta_i(\theta_{-i}|\theta_i) = \eta(\theta)/\eta_i(\theta_i)$ and
$\eta_i(\theta_i)$
is the marginal distribution of $\theta_i$.
\end{defi}

{\color{black}
In the literature of Bayesian games, $f_i:\Theta_i\to \mathcal{A}_i$ is called a {\em behavioural function}
and
consequently a Bayesian Nash equilibrium
 is also called a {\em behavioural function equilibrium}, see \cite{AgB06, Har67} and references therein.
Throughout this  paper, we will use both terminologies interchangeably for the equilibrium.
}

Note that there are a couple of alternative formulations for  (BNE).
If we let
\bgeqn
\rho_i(a_i,f_{-i},\theta_i) :=
%\min_{\eta_i(\theta_{-i}|\theta_i)\in \mathscr{Y}_{-i}}
\int_{\theta_{-i}\in \Theta_{-i}} u_i(a_i, f_{-i}(\theta_{-i}),\theta_i,\theta_{-i})d\eta_i(\theta_{-i}|\theta_i),
\label{eq:rho-i}
\edeqn
%{\color{blue} What is the definition of $\eta_i(\theta_{-i}|\theta_i)$????}
then we can
%rewrite
reformulate (BNE) as
\bgeqn
\inmat{(NE)}
  \quad f_i(\theta_i) \in\arg \max_{a_i\in \mathcal{A}_i} \rho_i(a_i,f_{-i},\theta_i), \;
  %\; \inmat{for}\;
  \forall \theta_i\in \Theta_i\, \inmat{and} \;  i\in N,
\label{eq:NE-max}
\edeqn
or
equivalently
\bgeqn
\rho_i(f_i(\theta_i),f_{-i},\theta_i) \geq \rho_i(g_i(\theta_i),f_{-i},\theta_i), \forall  \theta_i \in \Theta_i,
\label{eq:rho_i-BNE}
\edeqn
for every $g_i \in \mathcal{F}_i$, $i\in N$.
Consequently we may investigate existence of behavioural function equilibrium
of
(BNE)  by looking into
(NE).
For each $f\in \mathcal{F}$, define
\bgeqn
\Psi(f) :=\left\{(y_1(\cdot),\cdots, y_n(\cdot))\in \mathcal{F}:
y_i(\theta_i) \in\arg \max_{a_i\in \mathcal{A}_i} \rho_i(a_i,f_{-i},\theta_i), \forall \theta_i\in \Theta_i, i \in N\right\}.
\label{eq:Psi-f}
\edeqn
A sufficient condition for the well-definedness of $\Psi(f) $ is
compactness of $\mathcal{A}_i$ as well as continuity of $\rho_i$ in $a_i$ for $i\in N$.
On the other hand, if $u_i$ is concave and continuously differentiable w.r.t. $a_i$ for $i\in N$,
then $f$ is a behavioural function equilibrium if and only if
% (BNE) can be written as
it satisfies the following  variational inequality
\bgeqn
\label{eq:BNE-VIP}
%\inmat{(BNE)}  \quad
%f_i(\theta_i) \in\arg \max_{a_i\in A_i}
0\in \int_{\theta_{-i}\in \Theta_{-i}}
\nabla_{a_i} u_i(f_i(\theta_i), f_{-i}(\theta_{-i}),\theta_i,\theta_{-i})d\eta_i
(\theta_{-i}|\theta_i) + \mathscr{N}_{\mathcal{A}_i}(f_i(\theta_i)), \forall \theta_i\in \Theta_i, i\in N,
\edeqn
where $\mathscr{N}_X(x)$ denotes the normal cone of $X$ at point $x\in X$.
In what follows, we
%t might be helpful to
make a few comments on the definition of behavioural function equilibria and alternative formulations.

\begin{itemize}

\item[1.]
%First,
We require (\ref{eq:BNE}) to hold for every $\theta_i\in \Theta_i$, $i\in N$.
This differs from the Bayesian equilibrium model recently considered by
 Ewerhart \cite{EC14} and Ui \cite{Ui16} who require (\ref{eq:BNE}) to hold for almost every $\theta_i$ rather than every $\theta_i$ which means
that (\ref{eq:BNE}) may fail at a subset $\Theta_i^0$ of $\Theta_i$ with $\eta_i(\Theta_i^0)=0$.
% in his model.
 %Ewerhart \cite{EC14} calls the latter
 A behavioural function equilibrium defined in the ``almost sure'' sense is called a
 {\em pure strategy Nash equilibrium} (PSNE).
 %this kind of equilibrium
The difference will have a significant impact
%when we come towards
on conditions for existence and uniqueness of equilibria.  We will come back to this in Sections 3 and 4.
From the definition, we can see that a Bayesian behavioural function equilibrium is a pure Nash equilibrium but
not vice versa.
%it is not necessarily true conversely.
%From economic point of view, the former might be more interesting as
%there does not seem to
%{\color{blue}\textbf{Shaoyan, I think we may take a note to work on a model in a separate paper
%where (\ref{eq:BNE}) holds with some specified probability.
%I think it will be very interesting but have no idea whether it will be doable at this stage.
%}}
Note that
Meirowitz \cite{Mei03} does not make it clear on this but we can deduce from context of his paper that
his model also requires (\ref{eq:BNE}) to hold for every $\theta_i$.
% {\color{blue}
% Ewerhart \cite{EC14} calls the behavioural function equilibrium satisfying
%   (\ref{eq:BNE}) the
%    Bayesian Nash Equilibrium if
% (\ref{eq:BNE}) holds for every $\theta_i\in \Theta_i$, $i\in N$ and call that as
% pure strategy Nash Equilibrium if  (\ref{eq:BNE}) holds for almost every $\theta_i$.
%  }

\item[2.]
% Second,
 We implicitly
assume that maximum is attainable in each player's maximization problem (\ref{eq:BNE}).
This is guaranteed when $\mathcal{A}_i$ is compact and the expected utility function of each player is lower semi-continuous
w.r.t. its action variable.
%{\color{blue}
It is possible to replace the compactness condition
with  inf-compactness of the utility functions but we don't want
the additional technicality to  distract our focus on the key ideas.
%}

\item[3.]
%Third, a
An individual player may have multiple
global optimal solutions, denoted by $A^*_i(\theta_i)$, for some type values $\theta_i$,
in that case, $f_i(\cdot)$ is understood as a measurable
selection in the sense of Aumman \cite{Au65} from the set-valued mapping
$A^*_i:\Theta_i\to \mathcal{A}_i$.
% Fourth,
Moreover, we implicitly
assume
%s
that $u_i(a_i, f_{-i}(\theta_{-i}),\theta_i,\theta_{-i})$
is integrable with respect to $\eta_i(\theta_{-i}|\theta_i)$ over $\Theta_{-i}$.
A particularly interesting case is that $f_i$ is continuous on $\Theta_i$. We will focus on the case later on.

%%{\color{blue}
%Obviously $f_i:\Theta_i \to \mathcal{A}_i$ defines a function which is known as
%{\em  behavioural function} and hence a BNE is an $n$-tuple of
% Bayesian behavioural functions.
% %}
%%{\color{blue}
%%It is well known that $\eta(\theta_{-i}|\theta_i) = \frac{\eta(\theta)}{\eta_i(\theta_i)}$, %where
%%$\eta_i(\theta_i)$ denotes the marginal distribution of $\theta_i$.
%%}
%%One of the main theoretical issues concerning BNE is existence  of equilibria.
%%In this section, we explore the issue.

\item[4.] The behavioural function  equilibria
are not necessarily continuous. Indeed,
in some practical applications,
there might be a reason for discontinuity rather
than continuity, i.e., due to radical change
of technology in power generation or marketing strategy of a new product.
Here we give an academic example
with $A^*_i(\theta_i)$
being multi-valued and
(BNE) has
multiple discontinuous
behavioural function equilibria.

\end{itemize}

\begin{ex}[Multiple  discontinuous behavioural function equilibria]\rm
\label{ex:multiple-BNE-equi}
Let $u_1(a,\theta)=a_1 a_2 \theta_1$ and $u_2(a,\theta)=a_1 a_2 \theta_2$.
Let $\mathcal{A}_1=\mathcal{A}_2=[0,10]$ and $\Theta_1=\Theta_2=[-1,1]$.
Assume that $\theta_1$ and $\theta_2$ are  uniformly distributed over $\Theta_1$ and $\Theta_2$,
and $\theta_1$ and $\theta_2$ are independent.
We can easily figure out a behavioural function
equilibrium
$(f_1, f_2)$
with
\bgeqn
\label{eq:ex-f-1}
f_1(\theta_1)=
\begin{cases}
0         &\inmat{for}\, \theta_1\in [-1,0),\\
[0,10]  &\inmat{for}\, \theta_1=0,\\
%0         &\inmat{for}\, \theta_1\in [-1,0),
10       &\inmat{for}\, \theta_1\in (0,1],
\end{cases}
\edeqn
and
\bgeqn
\label{eq:ex-f-2}
f_2(\theta_2)=
\begin{cases}
0         &\inmat{for}\,\theta_2\in [-1,0),\\
[0,10]  &\inmat{for}\,\theta_2=0,\\
%0         &\inmat{for}\,\theta_2\in [-1,0).
10       &\inmat{for}\, \theta_2\in (0,1].
\end{cases}
\edeqn
Another behavioural function equilibrium is
$(f_1(\theta_1), f_2(\theta_2))=(0,0)$ for almost every
$(\theta_1, \theta_2) \in \Theta_1 \times \Theta_2$.

To see this, it follows from the definition of behavioural function equilibrium,
$(f_1,f_2)$ is an equilibrium if and only if
\bgeq
\displaystyle f_1(\theta_1) \in \arg\max_{a_1 \in [0,10]} \int_{-1}^1 \frac{1}{2}a_1 f_2(\theta_2) \theta_1 d \theta_2,
\forall \theta_1\in [0,1]
\edeq
and
\bgeq
\displaystyle f_2(\theta_2) \in \arg\max_{a_2 \in [0,10]} \int_{-1}^1 \frac{1}{2}a_2 f_1(\theta_1) \theta_2 d \theta_1,
\forall \theta_2\in [0,1].
\edeq
Since
$$
\int_{-1}^1 \frac{1}{2}a_1 f_2(\theta_2) \theta_1 d \theta_2
= \frac{1}{2}a_1\theta_1 \int_{-1}^1  f_2(\theta_2) d \theta_2,
$$
and
if $ \int_{-1}^1  f_2(\theta_2) d \theta_2>0$,
then
$f_1(\theta_1)=0$ for $-1 \leq \theta_1<0$,
$f_1(\theta_1)\in [0,10]$  for $\theta_1=0$, and
$f_1(\theta_1)=10$ for $0<\theta_1\leq 1$.
%for $-1 \leq \theta_1<0$, $f_1(\theta_1)=0$,
%for $\theta_1=0$, $f_1(\theta_1)\in [0,10]$ and
%for $0<\theta_1\leq 1$, $f_1(\theta_1)=10$.
%}
Likewise, we can obtain $f_2(\theta_2)$ as defined in (\ref{eq:ex-f-2}).
If $\int_{-1}^1  f_2(\theta_2) d \theta_2=0$,
then we can verify that $(f_1(\theta_1), f_2(\theta_2))=(0,0)$ for almost every
$(\theta_1, \theta_2) \in \Theta_1 \times \Theta_2$.
Obviously,  (BNE) has multiple discontinuous
behavioural function equilibria.
Note that in this example, we can see by the definition of PSNE that  there are two PSNEs.
\end{ex}

%{\color{red}
%\begin{rem}\rm
%\label{r:montonic-BNE}
%In the case when $z_i=1$, i.e., the dimension of $\mathcal{A}_i$ is 1,
%and $u_i$ is strictly monotonically increasing in $a_i$
%for fixed $a_{-i}$ and $\theta$,
%$i\in N$,
%then $\rho_i(a_i,f_{-i},\theta_i)$ is also a  strictly monotonically increasing function of
% $a_i$ for fixed $f_{-i}$ and $\theta_i$.
%Consequently, we may deduce from (\ref{eq:rho_i-BNE}) that
%any equilibrium
%$(f_1^*,\ldots,f_n^*)$
%of (BNE) or (NE)
%must satisfy
%\bgeqn
%\label{eq: per-senario}
%f_i^*(\theta_i) \geq g_i(\theta_i),\,\forall\theta_i \in \Theta_i, g_i \in \mathcal{F}_i, i \in N.
%\edeqn
%In other words, we may identify a behavioural function equilibrium through  (\ref{eq: per-senario}).
%
%\end{rem}
%}
%

\section{Existence of continuous behavioural function equilibrium}

In this section, we discuss the case when
each player's behavioural function is unique and continuous.
The uniqueness and continuity mean
that
each player's response is stable
against variation of its type (the behavioural function does not jump at any point of its domain).
In particular, we
 investigate conditions under which the behavioural function  equilibria are equi-continuous.
 The equi-continuity means that
 the derivatives of the
 player's behavioural
 functions are
 uniformly bounded. This is
   a key condition that Meirowitz used in his
    existence theorem \cite{Mei03} and he commented the condition is unlikely to be satisfied or verified.
    From computational point of view,
    the continuity allows us to develop efficient numerical schemes
    for solving (BNE),  which will be our focus in Section 4.

To this end, we need the following technical results about stability of a parametric programming problem.
To ease the notation, we will use $\|\cdot\|$ to denote the Euclidean norm in a finite dimensional space and any norm in a Banach space throughout the paper.

\begin{lema}[{\color{black} Quantitative stability of optimal solutions in parametric programming}]
\label{l-stability}
Let
$Z$ be a Banach space  equipped with norm $\|\cdot\|$,
$\phi, \psi:\R^m \times Z \rightarrow \R$
be continuous functions  and
$X \subseteq \R^m$
be a compact
set.
Consider
the following
parametric minimization problems
\begin{eqnarray}
\begin{array}{ll}
\displaystyle
\min_x & \phi(x,z) \\
\inmat{s.t.}\, & x\in X,
\end{array}
\label{eq:min-phi}
\end{eqnarray}
and
\begin{eqnarray}
\begin{array}{ll}
\displaystyle
\min_x & \psi(x,z) \\
\inmat{s.t.}\, & x\in X,
\end{array}
\label{eq:min-psi}
\end{eqnarray}
where $z\in Z$ is a parameter.
 For  $z_1, z_2\in Z$, let $X^*(z_1)$
 and $\tilde{X}^*(z_2)$
denote the set of optimal solutions to (\ref{eq:min-phi}) and (\ref{eq:min-psi}) respectively
with parameters $z_1$ and $z_2$.
Then
\begin{itemize}

\item[(i)] for any $\epsilon>0$, there exists a constant
$\delta>0$ (depending on $\epsilon$) such that when
$
\sup_{x\in X}|\phi(x,z_1)-\psi(x,z_2)|\leq \delta,
$
 \begin{equation}
 \mathbb{D}({\tilde X}^*(z_2),X^*(z_1)) \leq \epsilon,
 \label{eq:est-sol}
 \end{equation}
where $\mathbb{D}(B_1,B_2):=\sup_{b_1 \in B_1} d(b_1, B_2)$
with $d(b_1, B_2)= \inf_{b_2 \in B_2} \|b_1-b_2\|$;

\item[(ii)] if, in addition,
there exist positive constants $\alpha$ and $\nu$ such that
\begin{equation}
\phi(x,z_1) \geq \min_{x\in X} \phi(x,z_1) +\alpha d(x,X^*(z_1))^\nu, \;
\forall x\in X, \label{eq:2-grwth}
\end{equation}
then
\bgeqn
\mathbb{D}(\tilde{X}^*(z_2),X^*(z_1))
\leq
\left(\frac{3}{\alpha}\sup_{x\in X}|\psi(x,z_2)-\phi(x,z_1)|\right)^{\frac{1}{\nu}};
\label{eq:X-ptb-2order}
\edeqn

\item[(iii)]
%{\color{blue}
%if, further, $\phi(\cdot,z)$ is strictly quasiconvex,
%%$\phi(x,\cdot)$ is uniformly Lipschitz continuous over a set ${\cal Z}\subset Z$, that is,
%and
%there exists a positive constant $L$ such that
%$$
%|\psi(x,z_2)-\phi(x,z_1)| \leq L\|z_1-z_2\|, \forall z_1,z_2\in { Z},
%$$
%and
%the growth condition
%(\ref{eq:2-grwth}) holds for all $z \in {\cal Z}$,
%%}
%%the growth condition (\ref{eq:2-grwth}) holds for all $z\in {\cal Z}$,
%then
%$X^*(z)$ is a singleton for $z\in { Z}$ (written $X^*(z)=\{x^*(z)\}$ and
%%it
%$x^*(z)$
%is equi-continuous
%%over
% on ${\cal Z}$, that is,
%%there exists a positive constant such that
%\bgeqn
%\|x^*(z_2)-x^*(z_1)\|
%\leq
%\left(\frac{3}{\alpha}L\|z_2 -z_1\|\right)^{\frac{1}{\nu}}.
%\label{eq:X-equi-continu}
%\edeqn
%}
%{\color{red}
if, further, $\phi(\cdot,z)$ and $\psi(\cdot,z)$ are strictly quasiconvex,
%$\phi(x,\cdot)$ is uniformly Lipschitz continuous over a set ${\cal Z}\subset Z$, that is,
and
there exists a positive constant $L$ such that
$$
|\psi(x,z_2)-\phi(x,z_1)| \leq L\|z_1-z_2\|, \forall z_1,z_2\in { Z},
$$
and
the growth condition
(\ref{eq:2-grwth}) holds for all $z \in Z$,
%}
%the growth condition (\ref{eq:2-grwth}) holds for all $z\in {\cal Z}$,
then
$X^*(z)$ and $\tilde{X}^*(z)$ are singleton for $z\in { Z}$ (written $X^*(z)=\{x^*(z)\}$ and $\tilde{X}^*(z)=\tilde{x}^*(z)$ )
and
%it
%$x^*(z)$
%is equi-continuous
%over
% on ${\cal Z}$, that is,
%there exists a positive constant such that
\bgeqn
\|\tilde{x}^*(z_2)-x^*(z_1)\|
\leq
\left(\frac{3}{\alpha}L\|z_2 -z_1\|\right)^{\frac{1}{\nu}}.
\label{eq:X-equi-continu}
\edeqn
If $\psi=\phi$,
% in Lemma \ref{l-stability},
then (\ref{eq:X-equi-continu}) reduces to  the equi-continuity of the solution mapping $x^*(\cdot)$.
%}
\end{itemize}
\end{lema}

\noindent{\bf Proof.}
%Parts (i) and (ii).
% The results  are
%minor extensions of \cite[Lemma
%3.2]{dcb00} which deals with the case when $X^*(z_1)$ is a singleton
%and are also similar to \cite[Theorem 7.64]{RoW98}.
%Here we provide a proof for completeness.
Part (iii) follows directly from Part (ii), so we only prove Parts (i) and (ii).

Part (i). Let $\epsilon$ be a fixed
small positive number and $\phi^*_1$ be the optimal value of (\ref{eq:min-phi}) with parameter
% $z_2$.
$z_1$.
Define
 \bgeqn
 \label{eq:function-R}
 R(\epsilon) :=\inf_{\{x\in X:  d(x,X^*(z_1))\geq \epsilon\}}\phi(x,z_1)-\phi^*_1.
 \edeqn
Then $R(\epsilon)>0$. Let $\delta:=R(\epsilon)/3$ and $z_2$ be
such that $\sup_{x\in X} |\psi(x,z_2)-\phi(x,z_1)|\leq \delta$. Then for
any $x\in X$ with $d(x,X^*(z_1))\geq \epsilon$ and
for
any fixed
$x^*\in X^*(z_1)$,
  $$
\psi(x,z_2)-\psi(x^*,z_2)\geq \phi(x,z_1)-\phi(x^*,z_1)-2\delta \geq
R(\epsilon)/3> 0,
  $$
which implies that $x$ is not an optimal solution to
(\ref{eq:min-phi}) with parameter $z_2$.
This is equivalent to $d(x,X^*(z_1))<
\epsilon$ for all $x\in \tilde{X}^*(z_2)$, that is,
$\mathbb{D}(\tilde{X}^*(z_2),X^*(z_1)) \leq \epsilon$.

Part (ii). Under condition (\ref{eq:2-grwth}), it is easy to derive
that $R(\epsilon)=\alpha\epsilon^\nu$.  Let
$$
\epsilon :=\left(\frac{3}{\alpha}\sup_{x\in
X}|\psi(x,z_2)-\phi(x,z_1)|\right)^{\frac{1}{\nu}}.
$$
From Part (i), we immediately arrive at (\ref{eq:X-ptb-2order}).
%Part (iii) follows immediately from Part (ii).
The proof is complete.
\hfill $\Box$

We follow the line of Meirowitz \cite{Mei03} to use Schauder's fixed point theorem
for proving existence of equilibria in (NE). To this end,
we recall some relevant basic
definitions and results in functional analysis.

A set in a topological space is called {\em relatively compact} if its closure is compact.
Let $W$ be a Banach space and $T:W\to W$ be an operator.
The operator $T$ is said to be {\em compact} if it is continuous and
maps bounded sets into relatively compact sets.
The following result characterizes
relative compactness of a set in
functional spaces.
By the well-known Arzela-Ascoli theorem,
a set $D\subset W$ is relatively compact if and only if
the functions in $D$ satisfy
 the following two conditions:
(a) uniform boundedness,
that is,
$$
\sup_{f\in D} \|f\|_{\infty} <\infty,
$$
and
(b) equi-continuity, i.e.,
for every $\epsilon>0$, there exists
a constant $\delta>0$ such that
$$
%\forall \epsilon>0, \exists \delta>0 \; \inmat{such that} \;
\sup_{f\in D}
%{\color{blue}
\|f(s')-f(s'')\| \leq \epsilon, \forall s', s'' \; \inmat{with}\;
%}
\|s'-s''\|<\delta.
$$
{\color{black}
%which means
%that
Thus, if $K$  is a nonempty convex subset of a Hausdorff topological vector space
$V$  and $T$  is a continuous mapping of $ K$  into itself such that $T(K)$  is contained in a compact subset of
$K$, then $ T$
 has a fixed point. The following theorem precisely addresses this.
}

\begin{thm}[Schauder's fixed point theorem, 1930] If $M$ is a nonempty, closed, bounded, convex subset of a Banach space and $T: M\to M$ is a compact operator, then $T$ has a fixed point.
\label{t-Shauder}
\end{thm}

We now return to discuss existence of continuous equilibria in (NE) and make the following assumption.

\begin{assu}
\label{A:BNE-Main}
Consider problems (BNE) and (NE).
%Assume:
%\begin{itemize}
%\item[(a)]
For $i\in N$, the following conditions hold.
(a) $u_i(a,\theta)$ is continuous over $\mathcal{A} \times \Theta$ and
for each $f_{-i}, \theta_i$, $\rho_i(\cdot,f_{-i},\theta_i)$ is strictly
quasi-concave on $\mathcal{A}_i$;
(b) for a.e. measurable set $S \subset \Theta_{-i}$, $\eta_i(S|\theta_i)$ is continuous in $\theta_i$;
%({\color{blue}What role does this condition play??????});
% and for each $f_{-i}, \theta_i$, $\rho_i(\cdot,f_{-i},\theta_i)$ is strictly concave on $\mathcal{A}_i$;
 %(d) the optimal value of the maximization problem, denoted by $\vt_i(f_{-i},\theta_i)$, is continuous.
%\item[(b)]
(c) there exist positive constants $\alpha$ and $\nu$ such that
\begin{equation}
-\rho_i(a_i',f_{-i},\theta_i) \geq
-v_i(f_{-i},\theta_i) +\alpha d(a_i',A_i^*(f_{-i},\theta_i))^\nu, \;
\forall a'_i \in \mathcal{A}_i,
\label{eq:2-grwtha}
\end{equation}
where $v_i(f_{-i},\theta_i):=\max_{a_i \in \mathcal{A}_i}\rho_i(a_i,f_{-i},\theta_i) $ and
$A_i^*(f_{-i},\theta_i)
=\arg\max_{a_i \in \mathcal{A}_i}\rho_i(a_i,f_{-i},\theta_i)$;
%denotes the set of optimal solutions to
%the maximization problem in (NE);
 %\item[(c)] %{\color{blue}$\rho_i(a_i,f_{-i},\theta_i)$ is uniformly continuous in $(f_{-i},\theta_i)$ over $F_{-i}\times \Theta_i$.}
 (d)
 %{\color{blue}{[I change the order of words]}}
 there exists a positive constant
  $\tau_i>0$
  such that, for any $a_i \in {\mathcal{A}_i}$,$f_{-i} \in \mathcal{F}_{-i}$ and $\theta_i$, $\theta'_i \in \Theta_i$,
 \bgeqn
 \label{eq:lip-rho}
|\rho_i(a_i,f_{-i},\theta_i)-\rho_i(a_i,f_{-i},\theta'_i)|\leq \tau_i \|\theta_i-\theta'_i \|;
 \edeqn
% {\color{red} houlder?????}
 and (e)
 {\color{black}
 $\mathcal{A}_i$
and $\Theta_i$
 are compact and convex.}
\end{assu}
%{\color{red}
%Here and later on, $\|a\|$ denotes the Euclidian norm of $a$ when $a$ is a vector in a  finite dimensional space.
%}

Assumption \ref{A:BNE-Main} (a) is used by Meirowitz \cite{Mei03}, see conditions 2 and 3 in \cite[Proposition 1]{Mei03}.
It might be possible to weaken the continuity of $u_i$ in $a_i$ to lower semi-continuous but this would
incur more delicate analysis.  Assumption \ref{A:BNE-Main} (b) and (e)
coincide with conditions 5 and 1 respectively in \cite[Proposition 1]{Mei03}.
Assumption \ref{A:BNE-Main} (c) is newly introduced here. It requires $-\rho_i(a_i',f_{-i},\theta_i)$
to satisfy some
%kind of
growth condition at $A_i^*(f_{-i},\theta_i)$.
In the case when $\gamma=2$, this assumption is known as the second order growth condition which is
widely used in stability analysis of parametric programming, see \cite{BoS00}.
A sufficient condition for the latter is that $u_i$ is strongly concave in $a_i$ uniformly w.r.t. other parameters,
see Proposition \ref{p:second-growth}.
Assumption \ref{A:BNE-Main} (d) is also newly introduced here and requires $\rho_i$ to be uniformly Lipchitz continuous in $\theta_i$.
This condition may be weakened to H\"older continuity and we assume Lipschitz continuity only for the simplicity of presentation.
Note that the condition is satisfied if $u_i$ is uniformly equi-Lipschitz continuous in $\theta_i$ and
the density function $h_i(\cdot|\theta_i)$ of $\eta_i(\cdot|\theta_i)$ is Lipschitz continuous over $\Theta_i$, see Proposition \ref{l:uni-lip}.

\begin{thm}[Existence of continuous behavioural function equilibria] Consider problem (NE).
%Assume:
%%\begin{itemize}
%%\item[(a)]
%(a) $u_i(a,\theta)$ is continuous over $A \times \Theta$,
%for a.e. measurable set $S \subset \Theta_{-i}$, $\eta_i(S|\theta_i)$ is continuous in $\theta_i$,
% and for each $f_{-i}, \theta_i$, $\rho_i(\cdot,f_{-i},\theta_i)$ is strictly concave on $\mathcal{A}_i$;
% %(d) the optimal value of the maximization problem, denoted by $\vt_i(f_{-i},\theta_i)$, is continuous.
%%\item[(b)]
%(b) there exist positive constants $\alpha$ and $\nu$ such that
%\begin{equation}
%-\rho_i(a_i',f_{-i},\theta_i) \geq
%-v_i(f_{-i},\theta_i) +\alpha d(a_i',A_i^*(f_{-i},\theta_i))^\nu, \;
%\forall a'_i \in \mathcal{A}_i,
%\label{eq:2-grwtha}
%\end{equation}
%where $v_i(f_{-i},\theta_i):=\max_{a_i \in \mathcal{A}_i}\rho_i(a_i,f_{-i},\theta_i) $ and $A_i^*(f_{-i},\theta_i)$ denotes the set of optimal solutions to
%the maximization problem in (NE);
% %\item[(c)] %{\color{blue}$\rho_i(a_i,f_{-i},\theta_i)$ is uniformly continuous in $(f_{-i},\theta_i)$ over $F_{-i}\times \Theta_i$.}
% (c) for any $f_{-i} \in F_{-i}$ and $\theta_i$, $\theta'_i \in \Theta_i$, there exists
% a positive constant
%  $\tau_i$
%  %>0$
%  such that
% \bgeqn
% \label{eq:lip-rho}
%|\rho_i(a_i,f_{-i},\theta_i)-\rho_i(a_i,f_{-i},\theta'_i)|\leq \tau_i|\theta_i-\theta'_i|,
% \edeqn
% and (d) $A$ is compact.
%% \end{itemize}
 Let Assumption \ref{A:BNE-Main} hold. Then (NE) has an equilibrium
 %{\color{blue}
 with the behavioural functions being equi-continuous.
 %}
\label{t-NE-existence}
\end{thm}

\noindent
\textbf{Proof.} We use Theorem \ref{t-Shauder} to prove the result.
%Following similar analysis in Steps 1 and 2 in \cite[Proposition 1]{Mei03}.
Let
%$F$ and
$\mathcal{C}$ be defined as in (\ref{eq:F-C}) and $\Psi$ be defined as in (\ref{eq:Psi-f}).
%We verify that
 Note that $\mathcal{C}$ is a non-empty closed, bounded and convex set of a Banach space
 equipped with the infinity norm. In what follows, we verify that $\Psi: \mathcal{C} \to \mathcal{C}$
is a compact operator.

%{\color{red}
Observe first that for each
%$\theta_i \in \Theta_i$ and
$f_{-i}\in \mathcal{F}_{-i}$,
Assumption \ref{A:BNE-Main} (a) and (b)
%imply
ensure that the objective function $\rho_i(a_i,f_{-i},\theta_i)$ is continuous in $a_i$ and
$\theta_i$,
and strictly quasi-concave in $a_i$ for each fixed $\theta_i$.
Together with (e),
we have  $A^*_i(f_{-i},\theta_i)$
being non-empty and  a singleton.
 By classical stability results (see e.g.
\cite[Theorem 4.2.1]{BGKKT82}),
%upper semi-continuous in $\theta_i$
%and
%hence it
$A^*_i(f_{-i},\theta_i)$ is continuous
in $\theta_i$, which means
%This means
for any $f_{-i} \in \mathcal{F}_{-i}$,
$A^*_i(f_{-i},\cdot) \in \mathcal{C}_i$.
Moreover, since $\rho_i(a_i,f_{-i},\theta_i)$
is continuous in $(f_{-i},\theta_i)$, using the same stability argument,
we deduce that
$A^*_i(f_{-i},\theta_i): \mathcal{C}_{-i} \times \Theta_i \to \mathcal{A}_i$
is continuous.

On the other hand, under Assumption \ref{A:BNE-Main} (c) and (d), it follows from Lemma \ref{l-stability} that
the optimal solution of each maximization problem in (NE) is
equi-continuous on $\Theta_i$, that is,
\bgeqn
\label{eq:BNE-unique-proof}
|A^*_i(f_{-i},\theta_i)-A^*_i(f_{-i},\theta'_i)|\leq
\left(\frac{3\tau}{\alpha}\|\theta_i -\theta'_i\|\right)^{\frac{1}{\nu}}, \forall \theta_i,\theta_i'\in \Theta_i, \forall f_{-i}\in {\cal F}_{-i},
\edeqn
where $\tau:=\max_{i\in N}\{\tau_i\}$.
Since $\Theta_i$ is compact,
%{\color{blue}
 for any small positive number $\delta$, there exists
%$\delta$-net
a finite number of points $\theta_i^1,\cdots,\theta_i^M\in \Theta_i$ such that
 for every $\theta_i\in \Theta_i$, there exists $k\in \{1,\cdots,K\}$
such that
 $\|\theta_i-\theta_i^k\|\leq \delta$. Moreover, by the continuity of $A^*_i(f_{-i},\theta_i)$,
we may %require
set $f_{-i}'$ to be sufficiently close to $f_{-i}$ with
$\|A^*_i(f_{-i}',\theta_i^k)-A^*_i(f_{-i},\theta_i^k)\|\leq \epsilon$ for $k=1,\cdots,K$
with $\epsilon$ being a sufficiently small number.
By exploiting  the equi-continuity of  $A^*_i(f_{-i},\cdot)$, we have
 \bgeq
%\sup_{\theta_i\in \Theta_i}
 \|A^*_i(f_{-i}',\theta_i)-A^*_i(f_{-i},\theta_i)\|
&\leq&  \|A^*_i(f_{-i}',\theta_i)-A^*_i(f_{-i}',\theta_i^k)\| +\|A^*_i(f_{-i}',\theta_i^k)-A^*_i(f_{-i},\theta_i^k)\|\\
&&+\|A^*_i(f_{-i},\theta_i^k)-A^*_i(f_{-i},\theta_i)\|
\\
&\leq& 2\left(\frac{3\tau}{\alpha}\delta\right)^{\frac{1}{\nu}}+
\epsilon
\edeq
 and hence
%by the finite covering theorem, we can easily show  by exploiting
%the continuity of $A^*(\cdot,\cdot)$ and the equi-continuity of  $A^*(f_{-i},\cdot)$ that for any small positive number $\epsilon$
$$
\sup_{\theta_i\in \Theta_i} |A^*_i(f_{-i}',\theta_i)-A^*_i(f_{-i},\theta_i)|\leq 2\epsilon
$$
%when we set
for $\delta\leq \frac{\alpha}{3\tau}\left(\frac{\epsilon}{2}\right)^\nu$.
%}
%when $f_{-i}'$ is sufficiently close to $f_{-i}$.
This implies that $A^*_i(f_{-i},\cdot): {\cal C}_{-i}\to {\cal C}_i$ is continuous for each $i\in N$
and hence
 $\Psi: \mathcal{C} \to \mathcal{C}$ is a continuous operator.
Together with the compactness of $\mathcal{A}$,
this shows that $\Psi$ is a compact operator.

By Theorem \ref{t-Shauder}, (NE) has an equilibrium.
 Moreover, it follows from (\ref{eq:BNE-unique-proof}) that
 the behavioural function equilibria are equi-continuous.
\hfill $\Box$

%{\color{red} Shaoyan: It might be helpful if we include a discussion on uniqueness of equilibrium.
%}

%{\color{blue}
%\textbf{Why do we introduce the following?}
%
%For convenience,  we define
%\bgeqn
%\label{eq:lip-set}
%\mathcal{C}_0:=\left\{f: \Theta \to A: |f_i(\theta_i)-f_i(\theta'_i)| \leq \left(\frac{3\tau}{\alpha}\|\theta_i -\theta'_i\|\right)^{\frac{1}{\nu}}\,\inmat{for}\, i \in N\right\},
%\edeqn
%where $\alpha$ and $\tau$ are given in Assumption \ref{A:BNE-Main}.
%%Theorem \ref{t-NE-existence}.
%It is obvious that  $\mathcal{C}_0$ is a compact set.
%Under Assumption \ref{A:BNE-Main},
%%conditions of
%%Theorem \ref{t-NE-existence},
%any equilibrium of (NE) is contained in $\mathcal{C}_0$.
%}

Note that the growth condition (\ref{eq:2-grwtha}) is only
a sufficient condition to ensure
%the
equi-continuity of the behavioural functions.
In some particular
cases,
% the
 equi-continuity condition
 %can be verified directly,
 may be derived
without such a condition, see for instances
%such as
%the
 rent-seeking contests
 % considered
 in \cite{Fey08}.
 We will come back to this
 later on.
%Conditions
%(\ref{eq:2-grwtha}) and (\ref{eq:lip-rho}) are essential in Theorem \ref{t-NE-existence}.
%In what follows, we investigate circumstances under which
The following proposition states
that in the case when $\nu=2$, that is, the growth is of  second order,
condition (\ref{eq:2-grwtha})
may be derived from
strong concavity of $u_i$ in $a_i$.

{\color{black}
\begin{prop}[Sufficient conditions for the growth condition]
\label{p:second-growth}
Suppose that for $i\in N$,
$u_i(a,\theta)$ is Lipschitz continuous over ${\cal A} \times \Theta$
and for each $f_{-i}$ and $\theta$,
$u_i(\cdot, f_{-i}(\theta_{-i}),\theta):
\mathcal{A}_i\to \R$ is
strongly concave
on $\mathcal{A}_i$,
i.e., there exists a positive constant $\sigma_i$ such that
\bgeqn
\label{eq:stong-cav-u_i-1}
u_i(t a_i'+(1-t) a_i, f_{-i}(\theta_{-i}),\theta_i,\theta_{-i})
&\geq&
t u_i(a_i', f_{-i}(\theta_{-i}),\theta_i,\theta_{-i})+(1-t)  u_i(a_i, f_{-i}(\theta_{-i}),\theta_i,\theta_{-i})
\nonumber \\
&&
+\frac{\sigma_i}{2} t (1-t) \|a_i'-a_i\|^2, \,\forall a_i, a_i'\in \mathcal{A}_i, t \in [0,1].
\edeqn
Suppose that ${\cal A}$ is a convex set.
Then
$-\rho_i(a_i,f_{-i},\theta_i)$ satisfies  the second order growth condition (\ref{eq:2-grwtha}) with $\nu=2$.
\end{prop}
}

\noindent
\textbf{Proof.}
Observe first that the strong concavity of
$u_i(\cdot, f_{-i}$ $(\theta_{-i}),\theta_i,\theta_{-i})$
entails the strong concavity of $\rho(\cdot,f_{-i},\theta_i)$.
This
 can be deduced
from (\ref{eq:stong-cav-u_i}) by integrating on both sides
of the inequality with $\eta_i(\theta_{-i}|\theta_i)$ over $\Theta_{-i}$, i.e.,
%strong concavity of $u_i$ in $a_i$, that is, (
%and moreover
\bgeqn
\rho_i(t a_i'+(1-t) a_i, f_{-i},\theta_i)
&\geq&
t \rho_i(a_i', f_{-i},\theta_i) +(1-t)  \rho_i(a_i, f_{-i},\theta_i)
\nonumber\\
&&+\frac{\sigma_i}{2} t (1-t) \|a_i'-a_i\|^2, \,\forall a_i, a_i'\in \mathcal{A}_i.
\edeqn
Moreover, by \cite[Theorem 23.1]{Roc97}, the concavity and Lipschitz continuity imply
directional differentiability of $\rho_i$ in $a_i$.
%we have
Subtracting both sides of the inequality by $\rho_i(a_i, f_{-i},\theta_i)$
and then  dividing by $t$ and driving $t$ to $0$, we obtain
\bgeqn
\label{eq:directional-convex}
\rho_i(a_i',f_{-i},\theta_i)-\rho_i(a_i,f_{-i},\theta_i)
\leq
 (\rho_i)'_{a_i}(a_i,f_{-i},\theta_i;a_i'-a_i)
-\frac{\sigma_i}{2}\|a_i'-a_i\|^2.
\edeqn
%Here we are using the fact that the directional differentiability of $u_i$ implies
%the directional differentiability of $\rho_i$ and this is guaranteed
%by Clarke's theory (see discussions at \cite[page 79]{Cla83})
% under the concavity of $\rho_i$ in $a_i$.
%and
%Moreover, the strong concavity in $a_i$ ensures that $A_i^*(f_{-i},\theta_i)$ is singleton.
On the other hand,
the strong concavity in $a_i$ ensures that $A_i^*(f_{-i},\theta_i)$ is singleton.
By the first order optimality condition of $\rho_i$ at $A_i^*(f_{-i},\theta_i)$,
$$
(\rho_i)'_{a_i}(A_i^*(f_{-i},\theta_i),f_{-i},\theta_i; a'_i-A_i^*(f_{-i},\theta_i)) \leq 0,\,\forall a_i' \in \mathcal{A}_i.
$$
%where we write $\N_{A}(a)$ for the normal cone of set $A$ at point $a$,
Combining the inequality
%we obtain from
(\ref{eq:directional-convex}), we obtain
\bgeq
\rho_i(a_i',f_{-i},\theta_i)-\rho_i(A_i^*(f_{-i},\theta_i),f_{-i},\theta_i)
\leq
%\nabla_{a_i} \rho_i(A_i^*(f_{-i},\theta_i),f_{-i},\theta_i)^T(a_i'-A_i^*(f_{-i},\theta_i))\\
%&& -\sigma_i\|a_i'-A_i^*(f_{-i},\theta_i)\|^2\\
%&\leq&
-\frac{\sigma_i}{2}\|a_i'-A_i^*(f_{-i},\theta_i)\|^2, \; \forall a_i'\in \mathcal{A}_i,
\edeq
which
indicates  the second order growth of $-\rho_i(\cdot,f_{-i},\theta_i)$
 at $A_i^*(f_{-i},\theta_i)$.
%The proof is complete.
\hfill $\Box$

In the case when $u_i$ is  continuously differentiable, condition (\ref{eq:stong-cav-u_i-1})
%$a_i$, that is,
is equivalent to %there exists
existence of a positive constant $\sigma_i$ such that
for any fixed $a_i\in \mathcal{A}_i$
\bgeqn
&&u_i(a_i', f_{-i}(\theta_{-i}),\theta_i,\theta_{-i})
-
u_i(a_i, f_{-i}(\theta_{-i}),\theta_i,\theta_{-i})\nonumber\\
&&\leq
\nabla_{a_i}u_i(a_i, f_{-i}(\theta_{-i}),\theta_i,\theta_{-i})
^T(a_i'-a_i) - \frac{\sigma_i}{2}\|a_i'-a_i\|^2, \forall a'_i\in \mathcal{A}_i.
\label{eq:stong-cav-u_i}
\edeqn

%}

Condition
(\ref{eq:lip-rho}) also plays a crucial role in Theorem \ref{t-NE-existence}.
The proposition below shows that the condition may be derived from
Lipchitz continuity of $u_i$ in $(a,\theta)$ over  $\mathcal{A} \times \Theta$ and
the density function of $\eta_i(\cdot|\theta_i)$ is Lipschitz continuous over $\Theta_i$.
The latter is slightly strengthened from
Assumption \ref{A:BNE-Main} (b)
%Meirowitz \cite{Mei03} who
which requires the density function to be continuous rather than Lipschitz continuous.
%\textbf{Professor Xu: I change ``Meirowitz \cite{Mei03}'' to  our  Assumption \ref{A:BNE-Main} (b).}

\begin{prop}[Sufficient conditions for
%uniform Lipschitz continuity condition
 the validity of (\ref{eq:lip-rho})]
\label{l:uni-lip}
Assume: (a) $u_i$ is Lipschitz continuous over $\mathcal{A} \times \Theta$ with modulus $\kappa_i$;
% for a.e. measurable set $A \subset \Theta_{-i}$, $\eta_i(A|\theta_i)$
(b) the density function $h_i(\cdot|\theta_i)$ of $\eta_i(\cdot|\theta_i)$ is Lipschitz continuous over $\Theta_i$ with modulus $\gamma_i$,
that is,
\bgeqn
\label{eq:lip-density}
|h_i(\theta_{-i}|\theta'_i)-h_i(\theta_{-i}|\theta''_i)| \leq \gamma_i
\|\theta'_i-\theta''_i\|, \forall \theta'_i,\theta_i'' \in \Theta_i, \forall \theta_{-i} \in \Theta_{-i},
\edeqn
for $i\in N$,
{\color{black}and (c) ${\cal A}$ and $\Theta$ are compact.}
Then the uniform Lipschitz continuity condition  (\ref{eq:lip-rho})
 holds.
\end{prop}

%{\color{red}
%if (3.23) changes to holder, any new?
%}

\noindent
\textbf{Proof.}
By the definition of $\rho_i$, we have
\bgeq
&&|\rho_i(a_i,f_{-i},\theta'_i)-\rho_i(a_i,f_{-i},\theta''_i)|\\
%&=&\left|\int_{\theta_{-i}\in \Theta_{-i}} u_i(a_i, f_{-i}(\theta_{-i}),\theta'_i,\theta_{-i})d\eta_i(\theta_{-i}|\theta'_i)
%-\int_{\theta_{-i}\in \Theta_{-i}} u_i(a_i, f_{-i}(\theta_{-i}),\theta''_i,\theta_{-i})d\eta_i(\theta_{-i}|\theta''_i)\right|\\
&&\leq
\int_{\theta_{-i}\in \Theta_{-i}} \left| u_i(a_i, f_{-i}(\theta_{-i}),\theta'_i,\theta_{-i})-
 u_i(a_i, f_{-i}(\theta_{-i}),\theta''_i,\theta_{-i})\right|d\eta_i(\theta_{-i}|\theta'_i)\\
%&&+\left|\int_{\theta_{-i}\in \Theta_{-i}} u_i(a_i, f_{-i}(\theta_{-i}),\theta''_i,\theta_{-i})d\eta_i(\theta_{-i}|\theta'_i)
%-\int_{\theta_{-i}\in \Theta_{-i}} u_i(a_i, f_{-i}(\theta_{-i}),\theta''_i,\theta_{-i})d\eta_i(\theta_{-i}|\theta''_i)\right|\\
%&\leq&
%\kappa_i|\theta'_i-\theta''_i|\int_{\theta_{-i}\in \Theta_{-i}} d\eta_i(\theta_{-i}|\theta'_i)
&& \;\; +
\int_{\theta_{-i}\in \Theta_{-i}} |u_i(a_i, f_{-i}(\theta_{-i}),\theta''_i,\theta_{-i})(h_i(\theta_{-i}|\theta_i')-h_i(\theta_{-i}|\theta''_i))|d \theta_{-i}\\
&&\leq
\kappa_i\|\theta'_i-\theta''_i\|+
\gamma_i\|\theta'_i-\theta''_i\|\int_{\theta_{-i}\in \Theta_{-i}} |u_i(a_i, f_{-i}(\theta_{-i}),\theta''_i,\theta_{-i})|d \theta_{-i}\\
&&\leq(\kappa_i+
\gamma_i \Delta_i)\|\theta'_i-\theta''_i\|,
\edeq
where we set
$
\Delta_i:=\max_{a\in \mathcal{A},\theta \in \Theta} |u_i(a,\theta)|\int_{\theta_{-i}\in \Theta_{-i}}d \theta_{-i}.
$
%Therefore,
This shows condition  (\ref{eq:lip-rho})
% holds for
is fulfilled with
$\tau_i :=(\kappa_i+\gamma_i \Delta_i)$.
\hfill $\Box$

By the continuity of behavioural functions,
the behavioural function equilibrium has an alternative characterization.

\begin{thm}[Equivalent formulation of the BNE model]
\label{p:BNE-inequality}
{\color{black} Let Assumption \ref{A:BNE-Main} hold.
%Suppose condition (a) in Theorem \ref{t-NE-existence} holds.
Then
$f$ is a continuous
behavioural function equilibrium of} (BNE)
if and only if
%$f \in \mathcal{C}$
%and
%satisfies
%one of  the following conditions is fulfilled:
{\color{black} it satisfies}
\bgeqn
\bbe_\eta[u_i(f_i(\theta_i),f_{-i}(\theta_{-i}),\theta)] \geq \bbe_\eta[u_i(g_i(\theta_i),f_{-i}(\theta_{-i}),\theta)],
\forall g_i \in \mathcal{C}_i,
%\forall
\; \inmat{for} \; i \in N,
\label{eq:exante-1}
\edeqn
or {\color{black} equivalently}
\bgeqn
f \in \arg\max_{g \in \mathcal{C}} \sum_{i=1}^n
\bbe_\eta[u_i(g_i(\theta_i),f_{-i}(\theta_{-i}),\theta)].
\label{eq:exante-2}
\edeqn
\end{thm}

\noindent
\textbf{Proof.}
Under
%condition (a)
{\color{black} Assumption \ref{A:BNE-Main},} we know from
%the proof of
Theorem \ref{t-NE-existence}
%{\color{blue} [This theorem is about NE not BNE, I can not see why we can deduce the existence of a pure behavioural equilibrium????]}
that
%if
every
%pure strategy
%Bayesian Nash
behavioural function equilibrium $f$ of (BNE)
% then it
is a continuous function on $\Theta$.
%Note that
Moreover
$$
\bbe_\eta[u_i(f_i(\theta_i),f_{-i}(\theta_{-i}),\theta)]
=\int_{\Theta_i} \left\{\int_{\Theta_{-i}} u_i(f_i(\theta_i),f_{-i}(\theta_{-i}),\theta) d\eta_i(\theta_{-i}|\theta_i)\right\} d \eta_i(\theta_i),
$$
where $\eta_i(\theta_i)$ is the marginal probability distribution of $\theta_i$.

{\color{black}
The ``if'' part. Let $f \in \mathcal{C}$ and $f$ satisfies (\ref{eq:exante-1}).
We show that $f$ is a behavioural function equilibrium of (BNE).
Assume for the sake of a contradiction that
% follows from the fact that when
$f$ is not an equilibrium of (BNE).
%then
Then,
} there exist some
%$1\leq i \leq n$
$i\in N$ and  $g_{i} \in \mathcal{C}_i$ such that for some $\overline{\theta}_{i} \in \Theta_{i}$
\bgeq
\int_{\theta_{-i}\in \Theta_{-i}} u_i(f_i(\overline{\theta}_i), f_{-i}(\theta_{-i}),\overline{\theta}_i,\theta_{-i})d\eta_i(\theta_{-i}|\overline{\theta}_i)
<\int_{\theta_{-i}\in \Theta_{-i}} u_i(g_i(\overline{\theta}_{i}), f_{-i}(\theta_{-i}),\overline{\theta}_{i},\theta_{-i})d\eta_i(\theta_{-i}|\overline{\theta}_{i}).
\edeq
Here the deviation $g_i$ is picked up from $\mathcal{C}_i$ because
every behavioural function of player $i$ at the equilibrium is continuous.
%Since
Together with
{\color{black} Assumption \ref{A:BNE-Main},}
%condition (a) of Theorem \ref{t-NE-existence} as assumed,
% holds,
%this
the inequality above
 implies that there exists a neighborhood $\mathcal{B}_{\overline{\theta}_i}$ of $\overline{\theta}_i$ such that
\bgeqn
&&\int_{\mathcal{B}_{\overline{\theta}_i}}\int_{\theta_{-i}\in \Theta_{-i}} u_i(f_i(\theta_i), f_{-i}(\theta_{-i}),\theta_i,\theta_{-i})d\eta_i(\theta_{-i}|\theta_i)d \eta_i(\theta_i)\nonumber\\
&& <\int_{\mathcal{B}_{\overline{\theta}_i}}\int_{\theta_{-i}\in \Theta_{-i}} u_i(g_i(\theta_{i}), f_{-i}(\theta_{-i}),\theta_{i},\theta_{-i})d\eta_i(\theta_{-i}|\theta_{i}) d \eta_i(\theta_i).
\label{eq:Prop4.3-proof}
\edeqn
Thus we can construct a continuous function $\tilde{g}$
such that $\tilde{g}_i(\theta_i)$ satisfies inequality (\ref{eq:Prop4.3-proof})
%=g_i(\theta_i)$
%satisfies
for $\theta_i \in \mathcal{B}_{\overline{\theta}_i}$
and
%otherwise
$\tilde{g}_i(\theta_i)=f_i(\theta_i)$ outside the neighborhood.
%Let $\overline{f}:=(f_1,\ldots,f_{i-1},\tilde{g}_i, f_{i+1},\ldots,f_n)$.
Then we have
%{\color{red}
\bgeq
\bbe_\eta[u_i(f_i(\theta_i),f_{-i}(\theta_{-i}),\theta)] < \bbe_\eta[u_i(\tilde{g}_i(\theta_i),f_{-i}(\theta_{-i}),\theta)],
\edeq
which %a contradiction
contradicts the fact that
{\color{black} $f$ satisfies (\ref{eq:exante-1}).}
% $f$ is an equilibrium of the (BNE).

{\color{black}
The ``only if'' part. Let $f$ be a behavioural function equilibrium of (BNE), we show that it satisfies
(\ref{eq:exante-1}).
This}
 is obvious in that
 for any $\theta_i \in \Theta_i$
\bgeq
\int_{\theta_{-i}\in \Theta_{-i}} u_i(f_i(\theta_i), f_{-i}(\theta_{-i}),\theta_i,\theta_{-i})d\eta_i(\theta_{-i}|\theta_i)
\geq \int_{\theta_{-i}\in \Theta_{-i}} u_i(g_i(\theta_i), f_{-i}(\theta_{-i}),\theta_i,\theta_{-i})d\eta_i(\theta_{-i}|\theta_i).
\edeq
for $i=1,\ldots,n$ and by
integrating w.r.t $\theta_i$ on both sides of the inequality, we obtain (\ref{eq:exante-1}).

We now turn to prove that
the equivalence between
(\ref{eq:exante-2}) and (\ref{eq:exante-1}).
%Suppose
Let $f\in \mathcal{C}$ satisfy (\ref{eq:exante-1}).
By summing up w.r.t. $i$ on both sides of  (\ref{eq:exante-1}),
we
%can
immediately obtain (\ref{eq:exante-2}).
On the other direction,
%suppose tha
let $f$ satisfy (\ref{eq:exante-2})
%it does not satisfy
but not (\ref{eq:exante-1}).
 %This means
 Then there exist  $i \in N$ and a continuous function $g_{i}$ such that
\bgeq
\bbe_\eta[u_i(f_i(\theta_i),f_{-i}(\theta_{-i}),\theta)] < \bbe_\eta[u_i(g_i(\theta_i),f_{-i}(\theta_{-i}),\theta)].
\edeq
Let $\tilde{f}:=(f_1,\ldots,f_{i-1},g_i,f_{i+1},\ldots,f_n)$.
Then
%we have
\bgeq
\sum_{i=1}^n
\bbe_\eta[u_i(f_i(\theta_i),f_{-i}(\theta_{-i}),\theta)]
<
\sum_{i=1}^n
\bbe_\eta[u_i(\tilde{f}_i(\theta_i),f_{-i}(\theta_{-i}),\theta)],
\edeq
which leads to a contradiction
to (\ref{eq:exante-2}) as desired.
%The proof is complete.
\hfill $\Box$

%Proposition
Theorem \ref{p:BNE-inequality} enables
us to
%represent
recast (BNE) as follows:
an $n$-tuple $f:=(f_1,\cdots,f_n)$ is a continuous
 behavioural function equilibrium if
\bgeqn
\inmat{(BNE$'$)} \quad
f_i \in \arg \max_{g_i\in \mathcal{C}_i}
%\bbe_\eta[u_i(f_i(\theta_i),f_{-i}(\theta_{-i}),\theta)] \geq
\bbe_\eta[u_i(g_i(\theta_i),f_{-i}(\theta_{-i}),\theta)]
%\forall g_i \in C_i, \forall
\; \inmat{for} \; i \in N,
\label{eq:BNE'}
\edeqn
or equivalently
\bgeqn
\inmat{(BNE$''$)} \quad
f \in \arg \max_{g \in \mathcal{C} }
%\bbe_\eta[u_i(f_i(\theta_i),f_{-i}(\theta_{-i}),\theta)] \geq
\sum_{i=1}^n \bbe_\eta[u_i(g_i(\theta_i),f_{-i}(\theta_{-i}),\theta)].
%,
%\forall g_i \in C_i, \forall
%\; \inmat{for} \; i \in N,
\label{eq:BNE'-opti-form}
\edeqn
{\color{black}
The reformulation is possible because
we are restricting behavioural function equilibria of (BNE)
to continuous functions over $\Theta$ without affecting the nature of
the problem under Assumption \ref{A:BNE-Main}.
This is one of the key reasons that motivates us to focus on
continuous behavioural function equilibria rather than general equilibria.

Note that we can easily find a counter example that the reformulation fails to work without continuity of behavioural
function equilibrium.
To see this, let us
revisit Example \ref{ex:multiple-BNE-equi}. In that context,
if $\int_{-1}^1 f_2(\theta_2) d \theta_2>0$  and $\int_{-1}^1 f_2(\theta_2) d \theta_2>0$,
condition (\ref{eq:exante-1}) can be written as
\bgeq
\left(\int_{-1}^{1} f_1(\theta_1) \theta_1 d \theta_1\right) \left(\int_{-1}^1 f_2(\theta_2) d \theta_2\right)
\geq
\left(\int_{-1}^{1} g_1(\theta_1) \theta_1 d \theta_1\right) \left(\int_{-1}^1 f_2(\theta_2) d \theta_2\right),
\edeq
and
\bgeq
\left(\int_{-1}^{1} f_2(\theta_2) \theta_2d \theta_2 \right) \left(\int_{-1}^1 f_1(\theta_1) d \theta_1\right)
\geq
\left(\int_{-1}^{1} g_2(\theta_2) \theta_2d \theta_2 \right) \left(\int_{-1}^1 f_1(\theta_1) d \theta_1\right),
\edeq
or equivalently
\bgeqn
\label{eq:fail-1}
\int_{-1}^{1} f_1(\theta_1) \theta_1 d \theta_1 \geq
\int_{-1}^{1} g_1(\theta_1) \theta_1 d \theta_1
\edeqn
and
\bgeqn
\label{eq:fail-2}
\int_{-1}^{1} f_2(\theta_2) \theta_2d \theta_2 \geq
\int_{-1}^{1} g_2(\theta_2) \theta_2d \theta_2
\edeqn
for any $(g_1,g_2)  \in \mathcal{F}_1 \times \mathcal{F}_2$,
where $\mathcal{F}_1$ and $\mathcal{F}_2$ are the set of
measurable functions mapping from
$ [-1,1]$ to $[0,10]$.
Let $f_1^*(\theta_1)=0$ for $\theta_1 \in [-1,0)$
and $f_1^*(\theta_1)=10$ for $\theta_1 \in [0,1]$ except at point
$\theta_1=0.8$ where $f_1^*(0.8) =8$.
It is easy to
%know that
see that $(f_1^*,f_2^*)$ satisfies
(\ref{eq:fail-1}) and (\ref{eq:fail-2})
% hold but it is not
but it is not an equilibrium of (BNE).
Indeed, we can revise the value of $f_i^*$ at a set of points with Lebesgue measure zero
without affecting its satisfaction to (\ref{eq:fail-1}) and (\ref{eq:fail-2}).

%
%{\color{red}
%In the case when $z_i=1$ and $u_i$ is
%strictly
% monotonically increasing in $a_i$ for $i\in N$,
%we follow from Theorem \ref{p:BNE-inequality} and Remark \ref{r:montonic-BNE} that
%a continuous behavioural function equilibrium satisfying (\ref{eq:exante-1}) if and only if
%it satisfies (\ref{eq: per-senario}).
%}

The
%significance
importance of formulation (BNE$'$) compared to (BNE) is that
each player's expected utility is defined as the expected value of
its utility w.r.t. the joint probability distribution
$\eta(\theta)$ of the vector of type parameters
$\theta$ rather than the conditional probability distributions $\eta_i(\theta_{-i}|\theta_i)$.
This brings substantial convenience when we
discuss approximate schemes for solving (BNE) in the next section.
Formulation (BNE$''$) allows us to look into the equilibrium problem from optimization
perspective. We will use both formulations interchangeably later on depending on
which one
%brings more convenience
is more convenient to use in a context.
%{\color{blue}
%Before concluding this section,
In what follows, we use (BNE$''$) to derive conditions for the uniqueness of equilibrium.

\begin{thm}[Uniqueness of equilibrium]
\label{t-unique-BNE}
Let Assumption \ref{A:BNE-Main} (b)-(e) hold.
Assume: {\color{black}(a) for $i\in N$, $u_i(a,\theta)$ is Lipschitz continuous over ${\cal A} \times \Theta$}
%\footnote{{\color{blue}
%It is well-known that a concave function is locally Lipschitz continuous in the relative interior its domain
%and globally Lipschitz continuous over any compact and convex set of its domain.
%Here by saying ``Lipschitz continuous'' we really mean
% ``global Lipschitz continuous'' and  omit the word ``global'' for the simplicity of exposition because
% throughout the paper, we will not use local
% Lipschitz continuity.}}
 and
concave in $a_i$;
(b) for any $f',f''\in \mathcal{C}$ with $f'\neq f''$,
\bgeq
\int_\Theta  \sum_{i=1}^n \left[(u_i)'_{a_i}(f'_i(\theta_i),f'_{-i}(\theta_{-i}), \theta; f''_i(\theta_i)-f'_i(\theta_i))+
(u_i)'_{a_i}(f''_i(\theta_i),f''_{-i}(\theta_{-i}), \theta; f'_i(\theta_i)-f''_i(\theta_i))
\right] \eta(d \theta)>0.
%;
%\nonumber\\
%[H(f'(\theta),\theta) - H(f''(\theta),\theta)]^T(f'(\theta)-f''(\theta))\eta(d\theta) <0,
%\label{eq:unique-condition-dd}
\edeq
%(c) {\color{black}${\cal A}$ is a convex set.}
Then (BNE) possesses a unique equilibrium.
%which satisfies
%$$
%0\in \int_\Theta H(f(\theta),\theta)\eta(d\theta) +{\cal G}(f),
%$$
%where
%$
%{\cal G}(f) := \mathscr{N}_{{\cal C}_1}(f_1)\times \cdots\times \mathscr{N}_{{\cal C}_N}(f_N).
%$
%%and ${\cal C}_i$ is the projection of ${\cal C}$ on ${\cal A}_i$.

\end{thm}

%{\color{red}{Liwei's Proof of Theorem \ref{t-unique-BNE}.}}
%{\color{blue}{
\noindent\textbf{Proof.} Note that condition (a) is strengthened from Assumption \ref{A:BNE-Main} (a) and
hence under the condition and the rest of conditions in
Assumption \ref{A:BNE-Main},
we know from
% existence of equilibria is ensured by
Theorem \ref{t-NE-existence}
that the (BNE)  has an equilibrium.
 % under Assumption \ref{A:BNE-Main}.
In what follows, we show the uniqueness of the equilibrium.
Suppose for the sake of a contradiction that there are two
%different
distinct behavioural function equilibria denoted by
$f$ and $\tilde f$. Then by condition (b),
\bgeqn
&&\int_\Theta  \sum_{i=1}^n \left[(u_i)'_{a_i}(f_i(\theta_i),f_{-i}(\theta_{-i}), \theta; \tilde{f}_i(\theta_i)-f_i(\theta_i))+
(u_i)'_{a_i}(\tilde{f}_i(\theta_i),\tilde{f}_{-i}(\theta_{-i}), \theta;f_i(\theta_i)- \tilde{f}_i(\theta_i))
\right] \eta(d \theta)\nonumber\\
&&>0.
\label{eq:H-mon-1}
\edeqn
%Moreover,
{\color{black} 
On the other hand, following a similar argument to that in the proof of Proposition \ref{p:second-growth}, we know that
both $u_i$ and  $\rho_i$ are directionally differentiable w.r.t. $a_i$. Moreover, since
$-u_i$ is Clarke regular (see \cite[Definition 2.3.4]{Cla83}), it follows from formula (4) in page 79 of Clarke \cite{Cla83}
that
$$
(\rho_i)'_{a_i}(f_i(\theta_i),f_{-i},\theta_i;  \tilde{f}_i(\theta_i)-f_i(\theta_i))
=\int_{\Theta_{-i}} (u_i)'_{a_i}(f_i(\theta_i),f_{-i}(\theta_{-i}), \theta;  \tilde{f}_i(\theta_i)-f_i(\theta_i)).
d\eta_i(\theta_{-i}|\theta_i).
$$
}
%which enable us to deduce}
Consequently, we have
\bgeq
&&\int_\Theta  \sum_{i=1}^n (u_i)'_{a_i}(f_i(\theta_i),f_{-i}(\theta_{-i}), \theta;  \tilde{f}_i(\theta_i)-f_i(\theta_i)) \eta(d \theta)\\
&&=  \sum_{i=1}^n \int_{\Theta_i} \left[\int_{\Theta_{-i}} (u_i)'_{a_i}(f_i(\theta_i),f_{-i}(\theta_{-i}), \theta;  \tilde{f}_i(\theta_i)-f_i(\theta_i))
d\eta_i(\theta_{-i}|\theta_i) \right] d\eta_i(\theta_i)\\
&&=\sum_{i=1}^n \int_{\Theta_i} (\rho_i)'_{a_i}(f_i(\theta_i),f_{-i},\theta_i;  \tilde{f}_i(\theta_i)-f_i(\theta_i))d\eta_i(\theta_i)
\leq 0,
\edeq
%In the last equality, we exchange
%Here we are using the fact that the directional differentiability of $u_i$ implies
%the directional differentiability of $\rho_i$ and this is guaranteed
%by Clarke's theory (see discussions at \cite[page 79]{Cla83})
% under the concavity of $\rho_i$ in $a_i$.
%where the second equality follows from \cite{Cla83} under condition (a), see discussions at
%page 79 in the book.
%%holds
% due to the fact that
%$G'(a,v)=\int_\Theta g'(a,\theta;v) \eta(d \theta)$ with $G(a)=\int_\Theta g(a,\theta) \eta(d \theta)$
%when $g$ is a Lipschitz continuous and concave in $a$ for each $\theta \in \Theta$.
%To see this, since $g$ is Lipschitz continuous and concave in $a$, then $-g$ is regular,
%which ensures that $(-G)'(a,v)=\int_\Theta (-g)'(a,\theta;v) \eta(d \theta)$
%(see Page 79 in book [Optimization and nonsmooth analysis]),
%then $G'(a,v)=\int_\Theta g'(a,\theta;v) \eta(d \theta)$.
 %And t
 %The
 %first
where the last inequality
 is derived from
the first order  optimality condition of $\rho_i$ at $f_i(\theta_i)$.
Likewise, we can utilize the first order optimality condition of $\rho_i$
at $\tilde{f}_i(\theta_i)$ to establish
\bgeq
\int_\Theta  \sum_{i=1}^n (u_i)'_{a_i}(\tilde{f}_i(\theta_i),\tilde{f}_{-i}(\theta_{-i}), \theta; f_i(\theta_i)- \tilde{f}_i(\theta_i)) \eta(d \theta)
\leq 0.
\edeq
Combining the two inequalities above, we obtain
\begin{equation*}
\int_\Theta  \sum_{i=1}^n \left[(u_i)'_{a_i}(f_i(\theta_i),f_{-i}(\theta_{-i}), \theta; \tilde{f}_i(\theta_i)-f_i(\theta_i))+
(u_i)'_{a_i}(\tilde{f}_i(\theta_i),\tilde{f}_{-i}(\theta_{-i}), \theta;f_i(\theta_i)- \tilde{f}_i(\theta_i))
\right] \eta(d \theta)\leq 0,
\end{equation*}
which is a contradiction
{\color{blue} to (\ref{eq:H-mon-1}).}
%condition (b).
\hfill $\Box$
%}

In the case when $u_i$ is continuously differentiable in $a_i$,
condition (b)
%may be fulfilled by a stronger but more verifiable condition as follows:
%written as
%{\color{red}
%Q: $\nabla u_i(f_i(\theta_i),f_{-i}(\theta_{-i}),\theta)$ is a continuous function w.r.t. $\theta$ and $\Theta$ is a compact set.
%Does the above condition hold trivially?
%}
is equivalent to
\bgeqn
\int_\Theta [H(f'(\theta),\theta) - H(f''(\theta),\theta)]^T(f'(\theta)-f''(\theta))\eta(d\theta) <0,
\label{eq:H-monotone}
\edeqn
for any $f',f''\in \mathcal{C}$,$f'\ne f''$,
where $H(f,\theta):=(\nabla_{a_i} u_i(f_i(\theta_i),f_{-i}(\theta_{-i}),\theta): i\in N)$.
%which
Inequality (\ref{eq:H-monotone}) means that $H(\cdot,\theta)$ is diagonally strictly monotone in $f$ over $\mathcal{C}$.
% in the sense that

%
%
%
%
%\noindent
%\textbf{Proof.} Existence of equilibrium is ensured by Theorem \ref{t-NE-existence} under Assumption \ref{A:BNE-Main}.
%It therefore suffices to show uniqueness. Under condition (a),
%$H(f(\theta),\theta)$ is well defined. Moreover,
%$f\in \mathcal{C}$ is an equilibrium if and only
%% we can write down the first order optimality condition at an equilibrium
%it satisfies the following variational inequality:
%\bgeqn
%%f(\cdot) \in \arg \max_{g \in \mathcal{C} }
%%\bbe_\eta[u_i(f_i(\theta_i),f_{-i}(\theta_{-i}),\theta)] \geq
%0\in \int_\Theta H(f(\theta),\theta)\eta(d\theta) +\N_{\mathcal{C}}(f),
%%\sum_{i=1}^n \bbe_\eta[u_i(g_i(\theta_i),f_{-i}(\theta_{-i}),\theta)].
%%,
%%\forall g_i \in C_i, \forall
%%\; \inmat{for} \; i \in N,
%\label{eq:BNE'-opti-form-H}
%\edeqn
%where $\N_{\mathcal{C}}(f):=\N_{\mathcal{C}_1}(f_1) \times \cdots \times \N_{\mathcal{C}_n}(f_n)$.
%On the other hand, if we define scalar product
%$$
%\langle f, g\rangle :=\int_\Theta f(\theta)^Tg(\theta) \eta(d\theta), \forall f,g\in \mathcal{C},
%$$
%then $\mathcal{C}$ is a Hilbert space equipped with $\mathscr{L}_\infty$-norm.
%Condition (\ref{eq:H-monotone}) means that $H(\cdot,\theta)$ is strictly monotone on  ${\cal C}$.
%By \cite[Lemma 12.2]{HaKoSc}, the variational inequality (\ref{eq:BNE'-opti-form-H}) has at most one solution.
%The conclusion follows. \hfill $\Box$
%

At this point, it might be helpful to comment on the differences between the existence and uniqueness
results established by Ui \cite{Ui16} and our results.
First, Ui
%derives
demonstrated  existence and uniqueness of behavioural function equilibria by
converting (BNE) into an infinite dimensional variational inequality problem
and showing that the latter has a unique solution when
$u_i$ is continuously differentiable in $a_i$ for $i\in N$, and
$H(\cdot,\theta)$ is strictly monotone and
satisfies some coerciveness condition, see \cite[Proposition 2]{Ui16}.
In other words, in his work, existence and uniqueness are established in one go.
Here we
%derive
show  existence and uniqueness separately and our proof of existence
is similar to Meirowitz's proof which does not
require continuous differentiability of $u_i$ in $a_i$ or strict concavity of
$u_i$ in $a_i$;
second, the behavioural functions at an equilibrium in \cite{Ui16} are not necessarily continuous,
 and at an equilibrium (\ref{eq:BNE})
 is required to
 hold for almost every $\theta_i$ rather than for
%surely rather than
every $\theta_i$,
 the latter allows Ui to establish an equivalence formulation
analogous to (\ref{eq:exante-1}) without restricting the behavioural functions to be
continuous functions.
We retain a proof for our uniqueness result
%despite
since it is
%similar to
derived under a
%slightly
weaker condition
than that in \cite[Proposition 2]{Ui16}
and
%as
we have a different meaning of uniqueness.
%{\color{red}
%Moreover, it is unclear how to solve the infinite dimensional variational inequality problem.
%\textbf{Shaoyan, we may come back to this after reading more references.}
%}
Third, it is possible to relax the compactness of $\Theta$ and
strengthen the condition on $H(f,\theta)$ by making it
integrably bounded, we leave interested readers to explore as it is not our main focus here.

}

%Now we verify all the required condition above through the following example.

%{\color{red} 1. It might be more interesting to derive the BNE directly than
To conclude this section, we use an example to explain the existence and uniqueness results established in this section.
%showing its existence through the established theoretical results.}
Looking back Example \ref{ex:multiple-BNE-equi}, we find that all conditions in Assumption \ref{A:BNE-Main}
are satisfied except the strict concavity condition. To amend this, we include a second order term in each of
the utility function to make them strictly concave. This motivates us to consider the following example.

%{\color{red}
\begin{ex}[Uniqueness of continuous behavioural function equilibrium]
\label{ex:unique-BNE-equi}
\rm
%Let
{\color{black}
Consider a two player Bayesian game with
utility functions
 $u_1(a,\theta)=a_1 a_2 \theta_1-a_1^2$ and $u_2(a,\theta)=a_1 a_2 \theta_2-a_2^2$, action spaces
 $\mathcal{A}_1=\mathcal{A}_2=[0,10]$ and type sets $\Theta_1=\Theta_2=[-1,1]$.
Assume $\theta_1$ and $\theta_2$ are independent and uniformly distributed over $\Theta_1$ and $\Theta_2$ respectively.}
%and $\theta_1$ and $\theta_2$ are independent.
Then
(BNE) has a unique equilibrium
%is
$(f_1^*, f_2^*)$,
where
\bgeq
%\label{eq:ex-f-1-1}
f_1^*(\theta_1)= 0 \; \inmat{for} \;
%\begin{cases}
%\frac{1}{4}\theta_1 \int_{-1}^1 f_2(\theta_2) d \theta_2       &\inmat{for}\,
\theta_1\in [-1,1],
%,\\
%0         &\inmat{for}\, \theta_1\in [-1,0],
%\end{cases}
\edeq
and
\bgeq
%\label{eq:ex-f-1-1}
f_2^*(\theta_2)= 0 \; \inmat{for} \;
%\begin{cases}
%\frac{1}{4}\theta_1 \int_{-1}^1 f_2(\theta_2) d \theta_2       &\inmat{for}\,
\theta_2\in [-1,1].
%,\\
%0         &\inmat{for}\, \theta_1\in [-1,0],
%\end{cases}
\edeq
%and
%\bgeqn
%\label{eq:ex-f-2-2}
%f_2(\theta_2)=
%\begin{cases}
%\frac{1}{4}\theta_2 \int_{-1}^1 f_1(\theta_1) d \theta_1       &\inmat{for}\, \theta_2\in (0,1],\\
%0         &\inmat{for}\,\theta_2\in [-1,0].
%\end{cases}
%\edeqn
%This is because
%
To see this, it follows from the definition of behavioural function equilibrium,
%BNE,
$(f_1^*,f_2^*)$
%is
%
is an equilibrium if
and only if
%the only 2-tuples  which
it
satisfies
\bgeqn
f_1(\theta_1) \in  \arg \max_{a_1 \in [0,10]} \int_{-1}^1 \frac{1}{2}a_1 f_2(\theta_2) \theta_1 d \theta_2-a_1^2,\; \forall \theta_1\in  [-1,1]
\label{eq:BNE-ex2-1}
\edeqn
and
\bgeqn
f_2(\theta_2) \in \arg \max_{a_2 \in [0,10]} \int_{-1}^1 \frac{1}{2}a_2 f_1(\theta_1) \theta_2 d \theta_1-a_2^2, \; \forall \theta_2\in  [-1,1].
\label{eq:BNE-ex2-2}
\edeqn
%but this
%The ``if'' part can be verified easily. The ``only if
It is easy to verify that $(f^*_1, f^*_2)$ satisfies the above two conditions.
To see that this is the only solution, we note that since $a_i$ is restricted to take values in $[0,10]$,
 $\int_{-1}^1 f_i(\theta_i)  d \theta_i\geq 0$ for $i=1,2$. Thus
 from (\ref{eq:BNE-ex2-1}), $f_1(\theta_1)=0$ for $\theta_1\in [-1,0]$.
 Likewise from (\ref{eq:BNE-ex2-1}), $f_2(\theta_2)=0$ for $\theta_2\in [-1,0]$.
Moreover, for $\theta_1 \in[0,1]$, let
 $\alpha:=\int_{0}^1 f_2(\theta_2)  d \theta_2.$
If $\alpha>0$, then the optimal solution from (\ref{eq:BNE-ex2-1}) is
$$
f_1(\theta_1) = \frac{1}{4}\alpha\theta_1, \; \inmat{for}\; \theta_1\in [0,1].
$$
Substituting this to (\ref{eq:BNE-ex2-2}), we obtain
$$
f_2(\theta_2) = \frac{1}{32}\alpha\theta_2, \; \inmat{for}\; \theta_2\in [0,1].
$$
Substituting $f_2(\theta_2) $ back to (\ref{eq:BNE-ex2-1}), we obtain
$$
f_1(\theta_1) = \frac{1}{256}\alpha\theta_1, \; \inmat{for}\; \theta_1\in [0,1].
$$
Continuing the process, we deduce that $f_1(\theta_1)=f_2(\theta_2)=0$ for $\theta_1,\theta_2\in [0,1]$ in order for them to satisfy
conditions (\ref{eq:BNE-ex2-1}) and (\ref{eq:BNE-ex2-2}).

%for any $(\theta_1,\theta_2)\in [-1,1] \times [-1,1]$,
Note that
%To see that this is the only equilibrium, we resort to
the uniqueness can also be verified through Theorem \ref{t-unique-BNE}.
It is easy to calculate that
\bgeq
H(f,\theta) =
\left(
\begin{array}{c}
f_2(\theta_2)\theta_1 -2 f_1(\theta_1)\\
 f_1(\theta_1)\theta_2 -2 f_2(\theta_2))
\end{array}
\right)
=\left(
\begin{array}{cc}
-2  & \theta_1 \\
\theta_2 & -2
\end{array}
\right)
\left(
\begin{array}{c}
f_1(\theta_1)\\
f_2(\theta_2)
\end{array}
\right).
\edeq
Since the matrix at the right hand side of the equation is negative definite for every $(\theta_1,\theta_2)\in [-1,1]\times[-1,1]$,
then $H(\cdot,\theta)$ is diagonally strictly monotone.

\end{ex}

%Note that
%$$
%\int_{-1}^1 \frac{1}{2}a_1 f_2(\theta_2) \theta_1 d \theta_2-a_1^2
%= \frac{1}{2}a_1\theta_1 \int_{-1}^1  f_2(\theta_2) d \theta_2-a_1^2,
%$$
%and $ \int_{-1}^1  f_2(\theta_2) d \theta_2>0$,
%hence for $\theta_1>0$, $f_1(\theta_1)=\frac{1}{4}\theta_1 \int_{-1}^1 f_2(\theta_2) d \theta_2 $,
%and for $\theta_1\in [-1,0]$, $f_1(\theta_1)=0$.
%Likewise, we can obtain $f_2(\theta_2)$ as defined in (\ref{eq:ex-f-2-2}).
%
%Obviously, a continuous Bayesian Nash equilibria exists in this example.
%In what follows, we show that (NE) has a continuous equilibrium through Theorem \ref{t-NE-existence}.
%
%
%First, it is easy to verify that $u_i$ is strongly concave in $a_i$ for $i=1,2$.
%Now, we prove that $u_i$ is Lipschitz continuous.
%Let $a', a'' \in A$ and $\theta', \theta'' \in \Theta$.
%Then
%\bgeq
%|u_1(a',\theta')-u_1(a'',\theta'')|
%&=&|a'_1a'_2 \theta'_1-{a'}_1^2-a''_1a''_2 \theta''_1+{a''}_1^2|\\
%&\leq&|a'_1a'_2 \theta'_1-a''_1a''_2 \theta''_1|+|{a'}_1^2-{a''}_1^2|\\
%&\leq&|\theta'_1a'_1||a'_2-a''_2|+|\theta'_1a''_2||a'_1-a''_1|+|a''_1a''_2||\theta'_1-\theta''_1|+|{a'}_1+{a''}_1||{a'}_1-{a''}_1|.
%\edeq
%Since $A$ and $\Theta$ are compact,
%we have
%\bgeq
%|u_1(a',\theta')-u_1(a'',\theta'')|
%&\leq& 10|a'_2-a''_2|+10 |a'_1-a''_1|+10^2 |\theta'_1-\theta''_1|+20|{a'}_1-{a''}_1|\\
%&\leq& 10(2\sqrt{2}+\sqrt{2}) \|a'-a''\|+ \sqrt{2}10^2\|\theta'-\theta''\|.
%\edeq
%Likewise, we can obtain uniform Lipschitz continuity of $u_2$.

%%%%%%%%%%%%%%%%%%%%%
\section{
%Polyhedral decision rule and
Approximation schemes for (BNE$''$) }

In this section, we move on to discuss
approximation schemes for
% numerical methods for finding an equilibrium of
(BNE) which are ultimately aimed to provide some
numerical solution avenues for computing an
approximate
 behavioural function equilibrium.
%and we do so by developing some approximation schemes for its equivalent form
We do so via (BNE$''$) as our focus is on
those equilibria
%um
where
the behavioural functions are continuous.
Approximation is needed because  (BNE$''$) is an infinite dimensional stochastic equilibrium problem
which is
%very
in general
 difficult for us to obtain an exact equilibrium unless the problem has a very simple structure as
in Example \ref{ex:unique-BNE-equi}. To this end, we take two steps: (i) restrict
the space of behavioural functions to polynomial functions
and consequently  (BNE$''$) reduces to a finite dimensional stochastic equilibrium problem; (ii)
develop
%approximation
discretization  schemes for the stochastic equilibrium problem.
%We do so by adopting
%either for discretization
%when $\theta$ is continuously distributed.
% or
%for reducing the scenarios
%adopt
%
The approach
in step (i) is similar to the well-known polynomial decision rules
which have been recently developed
for solving
%{\color{red}
two-stage robust optimization problems \cite{BK12}
whereas
the approach in step (ii)
is well-known in stochastic programming but it is
% the first time
not often to be used
in
stochastic equilibrium problems except sample average approximation method \cite{XuZ13}.
In both approaches, we derive error bounds for the approximated equilibria.
%quantitative
%apply it to the
%.

\subsection{Polyhedral
%decision rule
behavioural function for (BNE$''$)}

To ease the exposition of technical results, we confine
ourself to the case
%Suppose
that $\Theta_i \subset \R$ and $\mathcal{A}_i =[a_i,b_i] \subset \R$
are compact intervals for $i=1,\ldots,n$
although the approximation schemes and technical results can be
%easily
extended to
the case when $\mathcal{A}_i$ and $\Theta_i$ are in multi-dimensional spaces.
Let $\xi(t):=(1,t,t^2,\ldots )$ be the sequence of monominals in $t \in \R$,
and denote by $\xi_d(t)$ the finite subsequence of the first $d+1$ elements of $\xi(t)$.
Thus, any polynomial of degree $d$ can be represented as $v^T \xi_d(t)$ for $v\in \R^{d+1}$.

%Denoted
Denote by $S_d^i$
%a
the set of polynomial functions
with
%of
the highest degree $d$:
\bgeq
S_d^i:=\left\{s: \Theta_i \rightarrow \mathcal{A}_i: \exists v \in \R^{d+1}\,\inmat{such that}\,s(t)=v^T \xi_d(t)\right\}
\edeq
for $i=1,\ldots,n$ and let $S_d:=(S_d^1,\ldots,S_d^n)$.
%Now we

We consider an approximation scheme for (BNE$''$)
%(\ref{eq:BNE'-opti-form})
by restricting
each player's behavioural functions
to $S_d^i$.
% rather than $\mathcal{C}_i$.
Consequently, we consider an $n$-tuple $f_d:=((f_d)_1,\cdots, (f_d)_n)$
%restrict the behavioural function  in $S_d$ and
%consider the following equilibrium problem:
such that
\bgeqn
\inmat{(BNE-app)} \quad f_d \in
\displaystyle\arg\max_{g_d \in S_d} \sum_{i=1}^n
\bbe_\eta[u_i((g_d)_i(\theta_i),(f_d)_{-i}(\theta_{-i}),\theta)].
\label{eq:poly-equi-1}
\edeqn
%{\color{red}
%or equivalently an $n$-tuple of $(d+1)$-dimensional vectors $v^*:=(v_1^*,\cdots,v_n^*)$ such that
%\bgeqn
%\inmat{(BNE-app')}  \quad v^* \in
%\displaystyle\arg\max_{v \in \R^{d+1}} \sum_{i=1}^n
%\bbe_\eta[u_i(v_i^T\xi_d(\theta_i),(v_i^*)^T\xi_d(\theta_i),\theta)].
%\label{eq:poly-equi-1'}
%\edeqn
%The significance of formulation (\ref{eq:poly-equi-1'}) is that
%it is a finite dimensional stochastic equilibrium problem
%which can be solved relatively more easily.
%{\bf Q: whether we need to remove the above to the place after (4.35).}
%}
{\color{black}
A
%The
significant
%ce
benefit  of formulation (\ref{eq:poly-equi-1})
%{eq:poly-equi-3})
is that
it is a {\em finite dimensional} stochastic equilibrium problem
which can be solved relatively more easily. }
To justify the approximation, we need to provide theoretical grounding which quantifies
the difference between an approximate equilibrium and its true counterpart.
%We will investigate the difference between problems  (\ref{eq:exante-2}) and  (\ref{eq:poly-equi-1}).
%Before this, we give a relationship between
We start by establishing a relationship between  $S_d$ and $\mathcal{C}$ in the following lemma.

\begin{lema}\label{l:dense}
The set $\{S_d\}$ is dense in $\mathcal{C}$ in the sense that for every $f\in \mathcal{C}$,
there exists a sequence $\{f_d\}\subset \{S_d\}$ such that
$\|f_d-f\|_\infty\to 0$ as $d$ tends to infinity.
\end{lema}

\noindent
\textbf{Proof.}
%Obviously,
Since polynomial functions are continuous, $S_d \subset \mathcal{C}$.
Without loss of generality, we
%suppose
assume that $\Theta_i =[0,1]$.
For any $f_i \in \mathcal{C}_i$, by the Weierstrass
theorem,
we can find a sequence of
Bernstein polynomials $B_d(\theta_i;f_i)$ of $f_i$, defined as
\bgeq
B_d(\theta_i;f_i):=
\sum_{j=0}^d f_i(j/d)
\binom{d}{j}
\theta_i^j (1-\theta_i)^{d-j}
\edeq
%satisfying that
such that
$
\|B_d(\theta_i;f_i)-f_i\|_{\infty} \to 0
$
as $d$ increases.
%It is easy to see that
Observe that
$$
B_d(\theta_i;f_i)\leq b_i
\sum_{j=0}^d
\binom{d}{j} \theta_i^j (1-\theta_i)^{d-j}=b_i,
$$
and likewise
$$
B_d(\theta_i;f_i)\geq a_i
\sum_{j=0}^d
\binom{d}{j} \theta_i^j (1-\theta_i)^{d-j}=a_i.
$$
This shows $B_d(\theta_i;f_i) \in S_d^i$
and the rest of conclusion is obvious.
%proof is competed.
\hfill $\Box$

%{\color{red}
%Question: $C$ is a a set of continuous function mapping from a compact $\Theta$
%to a compact $A$, in this case, whether  $S_d$ is dense in $C$.
%}
Based on Lemma \ref{l:dense},
we
are ready to show
%establish
 that
%the
any cluster point of
the sequence of
equilibria
%of Bayesian Nash equilibrium problem
obtained from solving  (\ref{eq:poly-equi-1})
is an equilibrium of (BNE$''$).
%Bayesian Nash equilibrium problem
% (\ref{eq:exante-2}).

\begin{thm}[Approximation of BNE by polynomial equilibria]
%Suppose
\label{t-BNE'-app}
Let $\{f_d\}$ be a sequence
of % an equilibrium of
approximate Bayesian behavioural function equilibria obtained from solving (BNE-app).
 %problem  (\ref{eq:poly-equi-1}).
%If $\|f_d-f\|_{\infty} \to 0$ as $d \to \infty$.
Then every cluster point of $\{f_d\}$
%$f$
 is
an equilibrium of Bayesian Nash equilibrium problem
 (\ref{eq:exante-2}).
\end{thm}

\noindent
\textbf{Proof.} Let $f$ be a cluster point and assume without loss of generality that
$\|f_d-f\|_{\infty} \to 0$. Since $f_d$ is an equilibrium of problem  (\ref{eq:poly-equi-1}),
then
%then for any $g_d \in S_d$,
\bgeqn
\label{eq:poly-equi-2}
\sum_{i=1}^n
\bbe_\eta[u_i((g_d)_i(\theta_i),(f_d)_{-i}(\theta_{-i}),\theta)]
\leq
\sum_{i=1}^n
\bbe_\eta[u_i(f_d(\theta),\theta)], \forall g_d \in S_d.
\edeqn
By Lemma \ref{l:dense},
polynomials are dense under the topology of
infinity norm in the space of continuous functions
on $\Theta$, which is  denoted by $\mathcal{C}$.
This means that
for any function $g\in \mathcal{C}$, there exists a sequence
%any cumulative point of sequences
of functions $\{g_d\} \subset \{S_d\}$ such that
%converging to it.
$\|g_d-g\|_\infty\to 0$ as $d\to \infty$.
%It follows from $\|f_d-f\|_{\infty} \to 0$ that $f \in C$.
Together with the continuity of $u_i$,
%by letting $d \to \infty$,
we obtain from (\ref{eq:poly-equi-2}) that for any $g \in {\cal C}$
\bgeqn
\label{eq:poly-equi-3}
\sum_{i=1}^n
\bbe_\eta[u_i(g_i(\theta_i),f_{-i}(\theta_{-i}),\theta)]
\leq
\sum_{i=1}^n
\bbe_\eta[u_i(f(\theta),\theta)],
\edeqn
which implies that
$f$ is an equilibrium of problem
 (\ref{eq:exante-2}).
%The proof is complete.
\hfill $\Box$

%{\color{blue}
Theorem \ref{t-BNE'-app} assumes the existence of polynomial equilibria in (BNE-app) for each fixed $d$.
In what follows, we investigate the existence.
%In order to solve problem
Let us rewrite (\ref{eq:poly-equi-1}) as:
%we rewrite it as follows:
\bgeqn
\inmat{(BNE-app$'$)}  \quad   V^* \in \displaystyle \arg\max_{V \in \mathcal{V}_d} \sum_{i=1}^n
\bbe_\eta[u_i(v_i^T \xi_d(\theta_i),(v^*_{-i})^T\xi_d(\theta_{-i}),\theta)],
\label{eq:poly-equi-3}
\edeqn
where $\xi_d(\theta_{-i}):=(\xi_d(\theta_j))_ {i \neq j \in N}$,
$V=(v_1,\ldots,v_{n}) \in \R^{(d+1)\times n}$,
$\mathcal{V}_d=(\mathcal{V}_d^1,\ldots,\mathcal{V}_d^n)$ with
$\mathcal{V}_d^i$ being
%a set in $\R^{d+1}$,
%
defined as
\bgeqn
\label{eq:V-d}
\mathcal{V}_d^i :=\left\{v_i \in \R^{d+1}: a_i \leq v_i^T \xi_d(\theta_i) \leq b_i,\, \forall \theta_i \in
\Theta_i\right\}.
\edeqn
%The significance of formulation (\ref{eq:poly-equi-3}) is that
%it is a {\em finite dimensional} stochastic equilibrium problem
%which can be solved relatively more easily.
The following lemma shows that  $\mathcal{V}_d^i$ is compact for $i \in N$.

%{\color{blue}
%Q: we should put a note that ${\cal V}^d_i$ defined in (4.35) is non-empty.
%}

\begin{lema}
\label{l:compact-V}
Let $\mathcal{V}_d^i$  be defined as in (\ref{eq:V-d}).
Then  $\mathcal{V}_d^i$ is a nonempty,  convex and compact set for
$i \in N$.
\end{lema}

\noindent
\textbf{Proof.}
Non-emptiness is obvious because we can always find a vector $v_i$ with the first component taking a value
between $a_i$ and $b_i$ and the other components being zero.
The convexity follows from the linear system of inequalities in $v_i$.
In what follows, we show compactness.

The closeness of $\mathcal{V}_d^i$ is obvious.
To see boundedness,
we select $d+1$ points $\theta_i^1,\theta_i^2,\ldots,\theta_i^{d+1} \in \Theta_i$
%and
with $\theta_i^j\neq \theta_i^k$ for $j\neq k$
and consider the following finite system of inequalities:
%constraints:
$$
a_i \leq v_i^T \xi_d(\theta_i^j) \leq b_i, j=1,\ldots,d+1.
$$
%which can be characterized as follows:
The system can be written in a matrix-vector form:
$$
a_i e \leq A({\theta_i}) v_i \leq b_i e,
$$
where $e$
%is vector
denotes
% in
the vector in $\R^{d+1}$ with
%all elements being $1$
unit components
and $A({\theta_i}) \in \R^{(d+1)\times(d+1)}$
%is called
%being
the Vandermonde matrix
 defined as
\bgeq
A({\theta_i}):=
\left[
\begin{array}{ccccc}
1 &\theta_i^1 &(\theta_i^1)^2 &\ldots &(\theta_i^1)^{d}\\
1 &\theta_i^2 &(\theta_i^2)^2 &\ldots &(\theta_i^2)^{d}\\
\vdots &\vdots &\vdots           & \ldots &\vdots\\
1 &\theta_i^{d+1} &(\theta_i^{d+1})^2 &\ldots &(\theta_i^{d+1})^{d}
\end{array}
\right].
\edeq
%By the property of the Vandermonde matrix,
%we know that
It is well known that $A(\theta_i)$ is nonsingular
and hence
%then
the set
$$
\tilde{\mathcal{V}}^i:=\{v_i \in \R^{d+1}: a_i e\leq  A(\theta_i)v_i \leq b_i e\}
$$
is bounded
{\color{black} because the null space defined by $A(\theta_i)v_i=0$ is $\{0\}$.}
%%Otherwise,
%To see the latter property,
%we assume for the sake of a contradiction that
%there exists a sequence $\{v_i ^M\}$ satisfying
%% \to \infty$ such that
% $$
% a_i e\leq  A(\theta_i) v_i^M \leq b_i e
% $$
% and
% $\|v_i^M\|\to \infty$ as $M$ tends to infinity.
%%by d
%Dividing both sides of the inequalities by $\|v_i^N\|$
%% on the above and letting
%and driving $M$
%%tend
% to $+\infty$,
%we have
% $$
%A(\theta_i) v_i^*=0
% $$
%with some $\|v_i^*\|=1$, which is a contradiction
%%with
%to  the nonsingularity of $A({\theta_i})$.
%Hence $\tilde{\mathcal{V}}^i$  is bounded.
Since $\mathcal{V}_d^i \subset \tilde{\mathcal{V}}^i$,
% and
%then
the boundedness of $\mathcal{V}_d^i$ is apparent.
%is obvious bounded.
%The proof is complete.
\hfill $\Box$

%Based on
By Lemma \ref{l:compact-V},
{\color{black}
we can find a convex and compact set  $\mathcal{U}_i \in \R^{d+1}$ such that $\mathcal{V}_d^i \in \mathcal{U}_i$
for $i \in N$. Consequently we may write
% then
\bgeqn
\label{eq:bounded-V-d}
\mathcal{V}_d^i=\left\{v_i \in \mathcal{U}_i: a_i \leq v_i^T \xi_d(\theta_i) \leq b_i,\, \forall \theta_i \in
\Theta_i\right\}
=\left\{v_i \in \mathcal{U}_i:v_i^T \xi_d(\theta_i)  \in S_d^i\right\}
\edeqn
and let
%Let
${\cal U}:=({\cal U}_1,\cdots,{\cal U}_n).$
}
%{\color{blue} Shaoyan: why don't we simply set $\mathcal{U}_i=\mathcal{V}_d^i$ and claim that $\mathcal{V}_d^i$ is a convex and compact set?
%In other words, we don't need to introduce $\mathcal{U}_i$ {\color{red}yes}.
% }
We are now ready to
 %establish the
 show existence of equilibria for (BNE-app$'$).
 % problem (\ref{eq:poly-equi-3}) in the following theorem.

\begin{thm}[Existence of polynomial equilibria in (BNE-app$'$)]
Suppose  that for $i\in N$,
 $u_i (a,\theta)$ is continuous over $\mathcal{A} \times \Theta$,
and $u_i(a_i,a_{-i},\theta)$ is concave in $a_i$.
Then problem (\ref{eq:poly-equi-3}) has an equilibrium.
\label{t-exist-poly-BNE}
\end{thm}

\noindent
\textbf{Proof.}
Let
$$
\phi(V,W)=\sum_{i=1}^n\bbe_{\eta}\left[u_i(v_i^T \xi_d(\theta_i),w_{-i}^T\xi_d(\theta_{-i}),\theta)\right],
$$
%{\color{red}
for $V=(v_1,\cdots,v_n) \in \R^{(d+1)\times n}$
and $W=(w_1,\cdots, w_n) \in  \R^{(d+1)\times n}$.
Since $u_i$ is continuous over $\mathcal{A} \times \Theta$ and concave in $a_i$,
then  $\phi(V,W)$ is continuous and concave in $V$ for each fixed $W$.
%and
%hence the e
Existence of an optimal solution to $\max_{V \in \mathcal{V}_d} \phi(V,W)$ is
ensured by
%the
compactness of $\mathcal{V}_d$.
To complete the proof,
we only need to show existence of $V^*$ such that
$$
V^* \in  \arg\max_{V \in \mathcal{V}_d} \phi(V,V^*).
$$
Denote by $\Upsilon(W)$
%be the
the set of optimal solutions to $\max_{V \in \mathcal{V}_d} \phi(V,W)$.
%for each fixed $W$.
Then $\Upsilon(W) \subset \mathcal{V}_d$.
By the concavity of $\phi(\cdot,W)$,
$\Upsilon(W)$ is a convex set.
Moreover,
it is easy to
%obtain
show  that $\Upsilon(W)$ is closed,
that is, for $W^k \to W^*$ and $U^k \in \Upsilon(W^k)$ with $U^k \to U^*$,
$U^* \in \Upsilon(W^*)$.
Furthermore,
it follows from \cite[Theorem 4.2.1]{AthLa06}
that $\Upsilon(W)$ is upper semi-continuous on $\mathcal{V}_d$.
By Kakutani's fixed point theorem \cite{Kak41},
there exists $V^*$ such that $V^* \in \Upsilon(V^*)$.
The proof is complete.
\hfill $\Box$

To see how (BNE-app$'$) works, we apply the approximation scheme to Example \ref{ex:unique-BNE-equi}.

%{\color{red}
\begin{ex}\rm
%Let us revisit Example \ref{ex:unique-BNE-equi} by
%using
%applying the polyhedral decision rule for
Let $d=1$, that is,
the behavioural functions are restricted to affine
functions
% the decision rule.
%Following the discussion in Example \ref{ex:unique-BNE-equi},
%we know that
%$f_1(\theta_1)=0$ for $\theta_1 \in [-1,0]$
%and
%$f_2(\theta_2)=0$ for $\theta_2 \in [-1,0]$.
%In what follows,
%we only consider the case that $\theta_1 \in (0,1]$ and $\theta_2 \in (0,1]$.
%$(f_1^*(\cdot),f_2^*(\cdot))=(
with $f_1(\theta_1) =v_0+v_1 \theta_1$
%+v_2\theta_1^2,
and $f_2(\theta_2) = w_0+w_1 \theta_2$.
%+ w_2 \theta_2^2)$
%is an equilibrium if
%and only if
%it
%satisfies
We need to find $(v^*_0, v_1^*, w_0^*,w_1^*)$ such that
\bgeqn
(v_0^*,v_1^*)
\in \arg \max_{(v_0,v_1) \in {\cal V}_1^1}
 \int_{-1}^{1}
 \left[\int_{-1}^1 \frac{1}{4} (v_0+v_1 \theta_1) (w_0^*+w_1^* \theta_2) \theta_1 d \theta_2-\frac{1}{2}(v_0+v_1 \theta_1)^2\right] d \theta_1
 %,
\label{eq:BNE-ex2-3}
\edeqn
and
\bgeqn
(w_0^*,w_1^*) \in \arg \max_{(w_0,w_1) \in {\cal V}_1^2}
\int_{-1}^{1}
\left[\int_{-1}^1 \frac{1}{4}(w_0+w_1 \theta_2)(v_0^*+v_1^* \theta_1)\theta_2 d \theta_1-\frac{1}{2}(w_0+w_1 \theta_2)^2
\right] d \theta_2,
\label{eq:BNE-ex2-4}
\edeqn
where
$$
{\cal V}_1^1:=\{(v_0,v_1): v_0+v_1 \theta_1\in [0,10], \theta_1 \in [-1,1]\}
%,
$$
and
$$
{\cal V}_1^2:=\{(w_0,w_1): w_0+w_1 \theta_2\in [0,10], \theta_2 \in [-1,1]\}.
$$
Problems (\ref{eq:BNE-ex2-3}) and (\ref{eq:BNE-ex2-4}) are constrained quadratic maximization problems.
Through some maneuvers, the equilibrium problem is down to solving
%we may simplify them to solve
%By calculation,
%we obtain that
%(\ref{eq:BNE-ex2-3}) and (\ref{eq:BNE-ex2-4}) can be simplified as follows:
\bgeqn
(v_0^*,v_1^*)
\in \arg \max_{(v_0,v_1) \in {\cal V}_1^1}
\frac{1}{3}v_1 w_0^*-v_0^2 -\frac{1}{3}v_1^2
%,
\label{eq:BNE-ex2-5}
\edeqn
and
\bgeqn
(w_0^*,w_1^*) \in \arg \max_{(w_0,w_1) \in {\cal V}_1^2}
\frac{1}{3}w_1 v_0^*-w_0^2 -\frac{1}{3}w_1^2.
\label{eq:BNE-ex2-6}
\edeqn
We can write down the KKT conditions for the two problems and solve the latter to
get an equilibrium. However, we opt for an easier way to identify an equilibrium.
For fixed $(w^1_0,w^1_1) \in {\cal V}_1^2$, we obtain from solving (\ref{eq:BNE-ex2-5}) that
$
v_0^1=v_1^1=\frac{w^1_0}{8}.
$
Substituting them to (\ref{eq:BNE-ex2-6}) and solve the latter, we obtain
$
w_0^2=w_1^2=\frac{v_0^1}{8}=\frac{w^1_0}{64}.
$
Continuing the process, we deduce that the only solution to (\ref{eq:BNE-ex2-5}) and (\ref{eq:BNE-ex2-6}) is
$(v_0^*,v_1^*,w_0^*,w_1^*)=(0,0,0,0)$ in that  $(0,0,0,0)$ satisfies (\ref{eq:BNE-ex2-5}) and (\ref{eq:BNE-ex2-6})
and the gradients of the two objective functions forms a mapping
%is
 $(-2v_0,\frac{1}{3}w_0-\frac{2}{3}v_1,-2w_0,\frac{1}{3}v_0-\frac{2}{3}w_1)$
which  is  strongly monotone
 (the Jacobian of the mapping is negative definite).

\end{ex}
%}

\subsection{
%Optimal
%approximation
%quantization method for
Discretization of (BNE-app$'$)}

%\textbf{Shaoyan: include some literature review on this, you may copy from your paper published in optimization and software}

%{\color{blue}
%The following paragraph has been re-organized.}
%In this part,
We now move on to discuss
% develop
discretization schemes for
%the stochastic equilibrium problem as discussed in the start of this section.
(BNE-app$'$) in the case when $\eta$ is continuously distributed.
%Here we select the
We consider the approach of
optimal quantization of probability measures due to Pflug and Pichler \cite{PfP11}
which
%is to find a
%find
identifies a discrete probability measure %which
approximating
%the true
% $P$
%probability measure optimally
$\eta$
under the Kantorovich metric.
Compared to the Monte Carlo methods and
Quasi-Monte Carlo methods,
this method has the highest approximation quality
with relatively fewer samples;
see comprehensive discussions by Pflug and Pichler \cite{PfP11}.
%}

%{\color{red}

%\textbf{Shaoyan: I copied the materials below from your Opti and Software paper.
%We need to edit as I don't think we need all of them.}

%\begin{defi}
%\label{d:kanto}
Let $\mathscr{L}$ denote
 the space of all Lipschitz continuous functions
$h:\Xi\rightarrow \R$ with Lipschitz constant no larger than $1$.
Let $P, Q \in \mathscr{P}(\Xi)$ be two probability measures.
Recall that
the Kantorovich metric (or distance) between $P$ and $Q$,
denoted by $\dd_K(P,Q)$, is defined by
\bgeq
\dd_K(P,Q):=\sup_{h \in \mathscr{L}}
          \left\{
          \int_{\Xi} h(\xi)P (d \xi)-
          \int_{\Xi} h(\xi)Q (d \xi)
          \right\}.
\edeq
%\end{defi}
%By Kantorovich-Rubinstein theorem \cite{KR58},
%\bgeq
%d_K(P,Q) =&\inf\left\{
%                      \int_{\Xi\times \Xi} \|\xi_1-\xi_2\| \pi(d \xi_1,d \xi_2):
%                      \begin{array}{l}
%                      \pi\;\inmat{is the joint distribution of }\;\xi_1\;\inmat{and}\;\xi_2\\
%                      \inmat{with marginals}\;P\;\inmat{and}\;Q,\inmat{respectively}
%                      \end{array}
%                     \right\}.
%\edeq
%The latter formulation is also known as Kantorovich formulation of
%%{\color{blue}
%Monge's transportation
%%}
%problem
%if we view  $P$ as goods spread over $\Xi$ to be relocated with new spread $Q$ over $\Xi$
%and $\|\xi_1-\xi_2\|$ as unit transportation cost \cite{Rac91}.
%Using the Kantorovich metric, we can
%quantify the distance between two sets of probability measures.
% Let $\CP,\CQ \subset \mathscr{P}$ be two sets of probability measures,
%we can define
%%{\color{blue}
%$$
%\mathbb{D}_K(\CP,\CQ)  := \sup_{P\in \CP}\inf_{Q\in \CQ} d_K(P,Q)
%$$
%%}
%which quantifies the deviation of $\CP$ from $\CQ$
%and
%$$
%\mathbb{H}_K(\CP,\CQ)  := \max\left\{\mathbb{D}_K(\CP,\CQ), \mathbb{D}_K(\CQ,\CP)\right\}
%$$
%that quantifies the distance between $\CP$ and $\CQ$.
Recall also that $\{P_N\}$ is said to converge to $P\in \mathscr{P}$ {\em weakly} if
   \bgeq
   \lim_{N\to\infty} \int_{\Xi}h(\xi)P_{N}(d\xi) = \int_{\Xi}h(\xi)P(d\xi),
\label{eq:lim-weak-covg}
   \edeq
for each bounded and continuous function $h:\Xi\rightarrow \R$.
An important property of the Kantorovich metric is that it metrizes weak convergence of probability measures
when the support set is bounded, that is, a sequence of probability measures
$\{P_N\}$ converges to $P$ weakly if and only if
%{\color{blue}
$\dd_K(P_N,P)\to 0$
as $N$ tends to infinity.
%\end{defi}
%}

%Note that
Let  $\Theta$ be a compact set and
%without loss of generality.
%let $\Theta:=[a_1,b_1]\times \cdots \times [a_n,b_n]$.
%Given a sequence
$\tilde{\Theta}^M:=\{\theta^k, k=1,\ldots,M\} \subset \Theta$ be a discrete subset of $\Theta$.
We can define the Voronoi
partition
$\{\Theta^1,\ldots,\Theta^M\}$ of $\Theta$,
where $\Theta^i$ are pairwise disjoint
with
\bgeq
\Theta^k
\subseteq
\left\{\theta \in \Theta: \|\theta-\theta^k\|=\min_{j=1,\ldots,M} \|\theta-\theta^j\|
\right\}.
\edeq
The possible optimal probability weights $p_k$
for minimizing $\dd_K\left(\eta,\sum_{k=1}^M p_k \delta_{\theta^k}\right)$ can then be found by
\bgeqn
\label{eq:optimal-p}
p=(p_1,\ldots,p_M)\;\inmat{with}\;p_k=\eta(\Theta^k).
\edeqn
Let
$
\eta^M(\cdot):=\sum_{k=1}^M p_k \delta_{\theta^k} (\cdot)
$
%where
with $p_k$
%is given
being defined as in (\ref{eq:optimal-p}).
%For  a set of probability measures $\A$ on $(\Xi, \mathscr{B})$,
%$\A$ is said to be {\em tight} if for any $\epsilon>0$, there exists a compact set $\Xi_\epsilon\subset \Xi$ such that
%$
%\inf_{P\in\A} P(\Xi_\epsilon)>1-\epsilon.
%$
%In the case when $\A$ is a singleton, it reduces to the tightness of a single probability measure.
%$\A$ is said to be {\em closed} (under the weak topology) if for any sequence $\{P_N\}\subset \A$
%with $P_N$ converging to $P$ weakly, we have $P\in\A$.
%$\A$ is said to be {\em weakly compact}
%if every sequence $\{P_N\}\subset \A$ contains a subsequence $\{P_{N'}\}$
%and $P\in \A$ such that $P_{N'}\rightarrow P$ weakly;
%%hx-edit
%see Billingsley \cite{Bil65}.
%
%
%By the well-known Prokhorov's theorem (see \cite{AthLa06}),
%a closed set $\A$ (under the weak topology) of probability measures is {\em compact} if it is tight.
%In particular, if $\Xi$ is a compact metric space,
%then the set of all probability measures on $(\Xi, \mathscr{B})$ is compact;
%see \cite[Theorem 1.12]{Pro56}.
%}
Following the definition of $\eta^M$ and the Kantorovich metric,
we have
\bgeq
\dd_K(\eta,\eta^M)= \int \min_{1 \leq k \leq M} d(\theta,\theta^k) d \eta(\theta)
=\sum_{k=1}^M \int_{\Theta^k} d(\theta,\theta^k) d \eta(\theta) \leq \beta_M,
\edeq
where $\beta_M$ is defined by
\bgeq
\label{eq:beta-M}
\beta_M:=\max_{\theta \in \Theta} \min_{1 \leq k \leq M} d(\theta,\theta^k) =\mathbb{H}(\tilde{\Theta}^M,\Theta).
\edeq
If  $\beta_M \to 0$ as $M \to \infty$,
then $\dd_K(\eta,\eta^M) \to 0$,
which implies that
$\eta^M$ converges to $\eta$ weakly.

Based on the discussions above, we may replace
$\eta$ with $\eta^M$
in
(BNE-app$'$)
and develop a discretization scheme for
the problem: find $V^M=(v_1^M,\cdots,v_n^M)$ such that
\bgeqn
\inmat{(\textbf{BNE-app$'$-dis})} \quad V^M \in \arg\max_{V \in \mathcal{V}_d^M} \sum_{i=1}^n
\bbe_{\eta^M}\left[u_i(v_i^T \xi_d(\theta_i),(v_{-i}^M)^T\xi_d(\theta_{-i}),\theta)\right],
\label{eq:poly-equi-saa}
\edeqn
where
\bgeqn
\label{eq:V-d-M}
(\mathcal{V}_d^M)^i :=\left\{v_i \in \mathcal{U}_i: a_i \leq v_i^T \xi_d(\theta_i) \leq b_i, \, \forall \theta_i \in \tilde{\Theta}_i^M\right\},
\edeqn
with $\tilde{\Theta}^M=\{\theta^j,j=1,\ldots,M\} \subset \Theta$ and
$\mathcal{U}_i$ being defined as in (\ref{eq:bounded-V-d}).
 In comparison with (BNE-app$'$), the (BNE-app$'$-dis) model only considers polynomial behavioural
functions defined over a discrete subset $\tilde{\Theta}^M$ of $\Theta$.
This might significantly enlarge the feasible set for $V$, that is, $\mathcal{V}_d^i \subset (\mathcal{V}_d^M)^i$.

Let $V^M$ be an equilibrium obtained from solving (BNE-app$'$-dis). In what follows, we investigate
convergence of $V^M$ as $M$ goes to infinity. To this end, we discuss convergence
of the feasible set of the maximization problem (\ref{eq:poly-equi-saa}) to that of (\ref{eq:poly-equi-3}), that is,
$(\mathcal{V}_d^M)^i$ to $\mathcal{V}_d^i$.
To ease the exposition, we introduce the following notation.

%
%
%For any two points $W=(w_1,\ldots,w_n) \in \R^{(d+1) \times n}$
%and $V=(v_1,\ldots,v_n) \in \R^{(d+1) \times n}$,
%%the distance between them is defined as follows:
%we let
%\bgeq
%\|W-V\|:=\sum_{i=1}^n \|w_i-v_i\|.
%\edeq
%to signify the distance between $W$ and $V$.
%%Before
%%we investigate the convergence of optimal solution of problem (\ref{eq:poly-equi-saa})
%%to that of problem (\ref{eq:poly-equi-3}),
%%we first establish the convergence of $\mathcal{V}_d^M \to \mathcal{V}_d$ as $M \to \infty$ in the following lemma.
%Let $g(v_i,\theta_i):=v_i^T \xi_d(\theta_i)$
%%Define
%$$
%\psi_i(v_i):= \sup_{\theta_i \in \Theta_i} g(v_i,\theta_i),
%\,\,
%\psi_i^M(v_i):= \sup_{j =1,\ldots,M} g(v_i,\theta_i^j),
%$$
%and
%$$
%\tilde{\psi}_i(v_i):= \sup_{\theta_i \in \Theta_i} -g(v_i,\theta_i),
%\,\,
%\tilde{\psi}_i^M(v_i):= \sup_{j =1,\ldots,M} -g(v_i,\theta_i^j).
%$$
%Define
%$$
%\mathcal{A}_i=\{v_i \in \mathcal{U}_i: \psi_i(v_i) \leq b_i \},\,\mathcal{A}_i^M=\{v_i \in \mathcal{U}_i: \psi_i^M(v_i) \leq b_i\}.
%$$
%and
%$$
%B_i=\{v_i \in \mathcal{U}_i: \tilde{\psi}_i(v_i) \leq -a_i\},\,B_i^M=\{v_i \in \mathcal{U}_i: \tilde{\psi}_i^M(v_i) \leq -a_i\}.
%$$
%Then we can write
%\bgeq
%\mathcal{V}^i_d&=&\mathcal{A}_i \cap B_i,\\
%(\mathcal{V}_d^M)^i&=&\mathcal{A}_i^M \cap B_i^M.
%%\mathcal{V}^i_d&=&\{v_i \in \mathcal{U}_i: \psi_i(v_i)\leq b_i, \tilde{\psi}_i(v_i)\leq -a_i\},\\
%%(\mathcal{V}_d^M)^i&=&\{v_i \in \mathcal{U}_i: \psi_i^M(v_i) \leq b_i, \tilde{\psi}_i^M(v_i) \leq -a_i\}.
%\edeq
%

For any two points $W=(w_1,\ldots,w_n) \in \R^{(d+1) \times n}$
and $V=(v_1,\ldots,v_n) \in \R^{(d+1) \times n}$,
%the distance between them is defined as follows:
we use
%\bgeq
$\|W-V\|:=\sum_{i=1}^n \|w_i-v_i\|$
%\edeq
to signify the distance between $W$ and $V$.
%Before
%we investigate the convergence of optimal solution of problem (\ref{eq:poly-equi-saa})
%to that of problem (\ref{eq:poly-equi-3}),
%we first establish the convergence of $\mathcal{V}_d^M \to \mathcal{V}_d$ as $M \to \infty$ in the following lemma.
Let $g(v_i,\theta_i):=v_i^T \xi_d(\theta_i)$,
%Define
$$
\psi_i(v_i):= \sup_{\theta_i \in \Theta_i} g(v_i,\theta_i),
\quad
\psi_i^M(v_i):= \sup_{j =1,\ldots,M} g(v_i,\theta_i^j),
$$
and
$$
\tilde{\psi}_i(v_i):= \sup_{\theta_i \in \Theta_i} -g(v_i,\theta_i),
\quad
\tilde{\psi}_i^M(v_i):= \sup_{j =1,\ldots,M} -g(v_i,\theta_i^j).
$$
Let
$$
\Omega_i :=\left\{v_i \in \mathcal{U}_i: \psi_i(v_i) \leq b_i \},\,\Omega_i^M=\{v_i \in \mathcal{U}_i: \psi_i^M(v_i) \leq b_i\right\}
$$
and
$$
\tilde{\Omega}_i :=\left\{v_i \in \mathcal{U}_i: \tilde{\psi}_i(v_i) \leq -a_i\},\,\tilde{\Omega}_i^M=\{v_i \in \mathcal{U}_i: \tilde{\psi}_i^M(v_i) \leq -a_i\right\}.
$$
Consequently we can write $(\mathcal{V}_d^M)^i$ and $\mathcal{V}_d^i$ respectively as
\bgeq
\mathcal{V}^i_d=\Omega_i \cap \tilde{\Omega}_i, \quad
(\mathcal{V}_d^M)^i=\Omega_i^M \cap \tilde{\Omega}_i^M.
%\mathcal{V}^i_d&=&\{v_i \in \mathcal{U}_i: \psi_i(v_i)\leq b_i, \tilde{\psi}_i(v_i)\leq -a_i\},\\
%(\mathcal{V}_d^M)^i&=&\{v_i \in \mathcal{U}_i: \psi_i^M(v_i) \leq b_i, \tilde{\psi}_i^M(v_i) \leq -a_i\}.
\edeq
%We
%$(\mathcal{V}_d^M)^i$ to $\mathcal{V}_d^i$.
In what follows, we estimate the difference between $\Omega_i^M \cap \tilde{\Omega}_i^M$ and $\Omega_i \cap \tilde{\Omega}_i$.

%{\color{red}
%
%For convenience,
%we define
%$$
%\mathcal{A}_i=\{v_i \in \mathcal{U}_i: \psi_i(v_i) \leq b_i \},\,\mathcal{A}_i^M=\{v_i \in \mathcal{U}_i: \psi_i^M(v_i) \leq b_i\}.
%$$
%and
%$$
%B_i=\{v_i \in \mathcal{U}_i: \tilde{\psi}_i(v_i) \leq -a_i\},\,B_i^M=\{v_i \in \mathcal{U}_i: \tilde{\psi}_i^M(v_i) \leq -a_i\}.
%$$
%
%}
%

\begin{prop}
\label{p:V-cover}
Let $\mathcal{V}_d^i$ and $(\mathcal{V}_d^M)^i$  be defined as in (\ref{eq:V-d}) and  (\ref{eq:V-d-M}).
Assume for $i\in N$ that there exist
%there exist
$v_i^* \in \mathcal{U}_i$ and a positive number $\alpha_i>0$ such that
\bgeqn
\psi_i(v_i^*)-b_i < -\alpha_i  \quad \inmat{and} \; \tilde{\psi}_i(v_i^*) +a_i<-\alpha_i.
\label{eq:Slater}
\edeqn
Then for $i\in N$
\bgeqn
\mathbb{H} ((\mathcal{V}_d^M)^i,\mathcal{V}_d^i) \leq \frac{\Delta_i L_i}{\alpha_i } \mathbb{H}(\tilde{\Theta}^M,\Theta),
\edeqn
where
$\Delta_i$ is the diameter of the set $\Omega_i \cap \tilde{\Omega}_i$,
and $L_i$ is the uniform Lipschitz modulus of $g(v_i,\theta_i)$ w.r.t. $\theta_i$ over ${\cal V}_d^i$.
\end{prop}

In order to prove the proposition, we need the following technical result.

\begin{lema}
\label{l:discrete-quantitative}
Let $T$ and $X$ be compact sets in some Banach spaces and $f(t,x): T \times X \to \R$
 be a continuous function.
 % mapping from $T \times X$ to $\R$ and
Let $X^K :=\{x^1,\cdots,x^K\}\subset X$ be
a discrete subset of $X$.
% and $\tilde{X}_1,\cdots,\tilde{X}_K$ be the Voronoi partition of $X$.
If $f$ is uniformly Lipschitz continuous in $x$ with modulus $L$, then
\bgeqn
\max_{t \in T}\left|\max_{x \in X} f(t,x) - \max_{k=1,\cdots,K}  f(t,x^k)\right| \leq
%L\max_{k=1,\cdots,K} \Delta_k \leq
L\mathbb{H}(X^K,X).
\edeqn
%where $\delta_k$ denotes the diameter of $\Theta^k$.
\end{lema}

\noindent
{\bf Proof.} Let $\tilde{X}_1,\cdots,\tilde{X}_K$ be the Voronoi partition of $X$. Then
$$
\max_{x \in X} f(t,x) = \max_{k=1,\cdots,K} \max_{x \in \tilde{X}_k} f(t,x).
$$
Let
\bgeqn
R(t):=\left|\max_{k=1,\cdots,K} \max_{x \in \tilde{X}_k} f(t,x) - \max_{k=1,\cdots,K}  f(t,x^k)\right|.
% \leq L\max_{k=1}^K\Delta_k
\label{eq:R(t)-1}
\edeqn
Note that for any two bounded sequences $\{a_k\}, \{b_k\}$, it is well known that
$$
\left| \sup_k a_k -\sup_k b_k\right| \leq \sup_k |a_k-b_k|.
$$
Thus, from (\ref{eq:R(t)-1}), we have
\bgeq
R(t) &\leq& \sup_{k=1,\cdots,K} \left|\max_{x \in \tilde{X}_k} f(t,x) - f(t,x^k)\right|=  \max_{k=1,\cdots,K} \left| f(t,x^k_*) - f(t,x^k)\right|,
\edeq
where $x^k_*$ denotes the maximizer of $f(t,\cdot)$ in the Voronoi cell  $\tilde{X}_k$.
Using the uniform Lipschtz continuity of $f$ in $x$, we have
\bgeq
R(t) %&\leq& \sup_{k=1}^K\left|\sup_{\theta \in \Theta^k} g(v,\theta) - g(v,\theta^k)\rght|\\
\leq  \sup_{k=1,\cdots,K}  L\|x_*^k-x^k\| = \sup_{k=1,\cdots,K}  Ld(x_*^k,X^K)
\leq L \mathbb{D}(X, X^K)= L\mathbb{H}(X, X^K),
%=L\beta_K,
\edeq
where the first equality
follows from the definition of the Voronoi partition that $x^k$ is the point from $X^K$ which is closest to $x_*^k$,
and the last equality is due to the fact that $X^K\subset X$.
%which is the conclusion.
\hfill $\Box$

\noindent
\textbf{Proof of Proposition \ref{p:V-cover}.}
Note that $g(v_i,\theta_i)$ is a polynomial function and
${\cal U}$
is a compact set.
Thus $g(v_i,\theta_i)$ is uniformly Lipschitz continuous w.r.t. $\theta_i$ over
${\cal U}_i$
with
%Lipschitz
modulus $L_i$.
By Lemma~\ref{l:discrete-quantitative},
we have
\bgeq
\sup_{v_i \in {\cal U}_i}
 |\psi_i(v_i)-\psi_i^M(v_i)| \leq L_i \mathbb{H}(\tilde{\Theta}^M,\Theta)
 %,
\edeq
and
\bgeq
\sup_{v_i \in {\cal U}_i}
 |\tilde{\psi}_i(v_i)-\tilde{\psi}_i^M(v_i)| \leq L_i \mathbb{H}(\tilde{\Theta}^M,\Theta).
\edeq
On the other hand, since $\psi_i$
and $\tilde{\psi}_i$ are convex functions, under Slater condition (\ref{eq:Slater}),
%and there exist $v_i^*$ and positive number $\alpha_i>0$ such that
%$\psi_i(v_i^*)-b_i < -\alpha_i$ and $\tilde{\psi}_i(v_i^*) +a_i<-\alpha_i$,
%hence there exists a positive constant $C_i^1$ such that
we have, by Robinson's error bound (\cite{Rob75}), that
for any $v_i \in \Omega_i^M \cap \tilde{\Omega}_i^M$,
\bgeq
d(v_i, \Omega_i\cap \tilde{\Omega}_i)
&\leq& C_i( \|(\psi_i(v_i)-b_i)_+\|+\|(\tilde{\psi}_i(v_i)+a_i)_+\|)\\
%&\leq& C_i^1 ( \|(\psi_i(v_i)-\psi_i^M(v_i))_+\|+\|(\psi_i^M(v_i)-b_i)_+\|
%+\|(\tilde{\psi}_i(v_i)-\tilde{\psi}_i^M(v_i))_+\|+\|(\tilde{\psi}_i^M(v_i)+a_i)_+\| )\\
&\leq&C_i ( \|(\psi_i(v_i)-\psi_i^M(v_i))_+\|
+\|(\tilde{\psi}_i(v_i)-\tilde{\psi}_i^M(v_i))_+\|)\\
&\leq&C_i ( \|\psi_i(v_i)-\psi_i^M(v_i)\|
+\|\tilde{\psi}_i(v_i)-\tilde{\psi}_i^M(v_i)\|)\\
&\leq& 2C_i L_i \mathbb{H}(\tilde{\Theta}^M,\Theta),
\edeq
where $C_i$ can be bounded by $\Delta_i /\alpha_i$ and we write $(t)_+$ for $\max(0,t)$ with $t\in \R$.
Thus
\bgeq
\mathbb{D}(\Omega_i^M \cap \tilde{\Omega}_i^M,\Omega_i\cap \tilde{\Omega}_i) \leq \frac{2\Delta_i L_i}{\alpha_i} L_i \mathbb{H}(\tilde{\Theta}^M,\Theta).
\edeq
%By Robinson’s Theorem on error bound for the convex system (see \cite{Rob75}),
%where $\Delta_i$ is the diameter of the set $\mathcal{A}_i \cap \mathcal{B}_i$.
Note that
%Due to
$\Omega_i\cap \tilde{\Omega}_i  \subset \Omega_i^M \cap \tilde{\Omega}_i^M$, that is,
%we know
%\bgeq
$\mathbb{D}(\Omega_i\cap \tilde{\Omega}_i,\Omega_i^M \cap \tilde{\Omega}_i^M) =0$.
%\edeq
Then
%Hence
\bgeq
\mathbb{H}(\Omega_i\cap \tilde{\Omega}_i,\Omega_i^M \cap \tilde{\Omega}_i^M)
=\mathbb{D}(\Omega_i^M \cap \tilde{\Omega}_i^M,\Omega_i\cap \tilde{\Omega}_i) \leq \frac{2\Delta_i L_i}{\alpha_i} \mathbb{H}(\tilde{\Theta}^M,\Theta).
\edeq
The proof is completed.
\hfill $\Box$

We are now ready to state the main convergence result of this section.
Let
$$
\phi(V,W) :=\sum_{i=1}^n\bbe_{\eta}[u_i(v_i^T \xi_d(\theta_i),w_{-i}^T\xi_d(\theta_{-i}),\theta)]
$$
and
$$
\phi^M(V,W) :=\sum_{i=1}^n \bbe_{\eta^M}[u_i(v_i^T \xi_d(\theta_i),w_{-i}^T\xi_d(\theta_{-i}),\theta)].
$$

%{\color{blue}
We are now ready to state the main result of this section which describes
convergence  of approximate polynomial behavioural function equilibria solved from (BNE-app$'$-dis).
%}

\begin{thm}[Convergence of approximate polynomial behavioural function equilibria
%stochastic approximation
]
\label{t-V^MtoV}
Let $\{V^M\}$ be a sequence of equilibria obtained from solving (BNE-app$'$-dis).
Assume for $i\in N$: (a)
{\color{black}
$u_i(a_i,a_{-i},\theta)$ is Lipschitz continuous
% w.r.t.
in  $(a_i,a_{-i},\theta)$ with modulus $\kappa_i$,
}
%{\color{blue}
%$u_i(a_i,a_{-i},\theta)$ is
%uniformly H\"older continuous w.r.t. $a_{-i}$, i.e., there exists positive constant
%$\kappa_i$
%such that
%$$
%|u_i(a_i,a_{-i}',\theta) - u_i(a_i,a_{-i}'',\theta)| \leq \kappa_i\|a_{-i}'-a_{-i}''\|, \forall a_i\in {\cal A}_i,
%a_{-i}',a_{-i}''\in {\cal A}_{-i},\theta\in\Theta,
%$$ }
(b)
{\color{black}conditions in Proposition \ref{p:V-cover} hold and }
$\mathbb{H}(\tilde{\Theta}^M,\Theta)\to 0$ as $M$ tends to infinity.
Then
\begin{itemize}
\item[(i)]  any cluster point $V^*$ of the sequence $\{V^M\}$
%Any cumulative point  of  equilibriums of approximate problem (\ref{eq:poly-equi-saa})
is an equilibrium of (BNE-app$'$);
%a of problem (\ref{eq:poly-equi-3}).
%Then
%\bgeq
%\limsup_{M \to \infty} V^M \subset \mathcal{V}^*.
%\edeq
%Furthermore,

\item[(ii)] {\color{black}
%Moreover
if, in addition,  $\phi(V,V^*)$ satisfies
%the second
some growth condition at $V^*$,
that is,
there exist constants $\alpha>0$ and $0<\nu<1$ such that
\bgeq
-\phi(W,V^*) \geq -\phi(V^*,V^*)+ \alpha\|W-V^*\|^\nu, \forall W \in \mathcal{V}_d,
\edeq
and
then there exists a positive constant $C$ such that
$$
\|V^M-V^*\| \leq C(\dd_K(\eta,\eta^M))^{\frac{1}{\nu}}
$$
when $M$ is sufficiently large;
%\textbf{Let's keep this for the time being in case we have some new idea}

\item[(iii)] if $\theta^1,\cdots,\theta^M$ are independent and identically distributed with
%id samples and
$p_k=\frac{1}{M}$ for $k=1,\cdots,M$,
then for any small positive number $\varepsilon$, there exist positive constants $C(\varepsilon)$ and
$\lambda(\varepsilon)$ such that
\bgeq
\prob\left(\|V^M-V^*\| \geq \varepsilon\left| \|V^*-V^M\|\leq \left(\frac{\alpha}{4\kappa(n-1) K}\right)^{\frac{1}{1-\nu}}\right.\right) \leq C(\varepsilon)e^{-\lambda(\varepsilon)M}
\edeq
for  $M$ sufficiently large,
where $\kappa :=\max_{i=1,\cdots,n}\kappa_i$ and
$K:=\sum_{i=1}^n \sup_{\theta_i \in \Theta_i}\|\xi_d(\theta_{i})\|$.
In the case when
\bgeq
\prob\left(\|V^*-V^M\|\leq \left(\frac{\alpha}{4\kappa(n-1) K}\right)^{\frac{1}{1-\nu}}\right)=1
%\label{eq:Prob:V_M-V*=1}
\edeq
for $M$ sufficiently large,
\bgeq
\prob\left(\|V^M-V^*\| \geq \varepsilon
%\left| \|V^*-V^M\|\leq \left(\frac{\alpha}{4\kappa K}\right)^{\frac{1}{\beta-\nu}}\right.
\right) \leq C(\varepsilon)e^{-\lambda(\varepsilon)M}.
\edeq
Here the probability measure ``Prob'' is understood as the product
  of
  the true (unknown) probability measure of $P$ over the measurable space
  $\Xi\times\Xi\times\dots$ with product Borel sigma-algebra $\mathscr{B}\times\mathscr{B}\times\dots$.
}
\end{itemize}

\end{thm}

\noindent
\textbf{Proof.}
%Let $V^M$ be a sequence of equilibriums of approximate problem (\ref{eq:poly-equi-saa}).
Assume without loss of generality (by taking a subsequence if necessary)
that $\{V^M\}$ converges to $V^*$.
Let us write
\bgeqn
\label{eq:uniform-conver}
\phi(V,V^*)-\phi^M(V,V^M)=\phi(V,V^*)-\phi^M(V,V^*)+\phi^M(V,V^*)-\phi^M(V,V^M).
\edeqn
We estimate the right hand side of (\ref{eq:uniform-conver}).
Since $u_i$ is Lipschitz continuous and  $v^T \xi_d(\theta)$ is a polynomial function,
%we have
$u_i(v_i^T \xi_d(\theta_i), v_{-i}^T \xi_d(\theta_{-i}),\theta)$
is uniformly Lipschitz continuous
w.r.t. $\theta$ for all
%over compact set
$v\in \mathcal{U}$.
We denote the Lipschitz modulus by $L_i$.
Thus %for any $V \in \mathcal{U}$,
\bgeqn
&&\sup_{V \in \mathcal{U}}
|\phi(V,V^*)-\phi^M(V,V^*)| \nonumber \\
&&\leq \sup_{V \in \mathcal{U}}\sum_{i=1}^n \left|\bbe_{\eta}[u_i(v_i^T \xi_d(\theta_i),(v_{-i}^*)^T\xi_d(\theta_{-i}),\theta)]
-\bbe_{\eta^M}[u_i(v_i^T \xi_d(\theta_i),(v_{-i}^*)^T\xi_d(\theta_{-i}),\theta)]\right|\nonumber\\
&&\leq
%\sum_{i=1}^M L_i
L\dd_K(\eta,\eta^M),
\label{eq:uniform-conver-2}
\edeqn
where $L=\sum_{i=1}^n L_i$.
The last inequality follows from the definition of the Kantorovich metric.
%Moreover, we have
%\bgeqn
%\label{eq:uniform-conver-2}
%\sup_{V\in \mathcal{U}}|\phi(V,V^*)-\phi^M(V,V^*)| \leq   L \dd_K(\eta,\eta^M),
%\edeqn
%with $L=\sum_{i=1}^n L_i$
Now we turn to estimate $\phi^M(V,V^*)-\phi^M(V,V^M)$.
%the second term of (\ref{eq:uniform-conver}).
\bgeq
&&|\phi^M(V,V^*)-\phi^M(V,V^M)|\\
&&\leq
\sum_{i=1}^n \left|\bbe_{\eta^M}\left[u_i(v_i^T \xi_d(\theta_i),(v_{-i}^*)^T\xi_d(\theta_{-i}),\theta)\right]
-\bbe_{\eta^M}\left[u_i(v_i^T \xi_d(\theta_i),(v_{-i}^M)^T\xi_d(\theta_{-i}),\theta)\right]\right|\\
&&\leq \sum_{i=1}^n \kappa_i  \bbe_{\eta^M}\left[|(v_{-i}^*-v_{-i}^M)^T\xi_d(\theta_{-i})|\right] \quad \inmat{(using Lipschitz
%H\"older
continuity of $u_i$)} \\
&&\leq \kappa (n-1) \sum_{i=1}^n  |(v_{i}^*-v_{i}^M)^T \bbe_{\eta^M}[\xi_d(\theta_{i})]|\\
&&\leq \kappa (n-1) K  \|V^*-V^M\|,
\edeq
where $\kappa :=\max_{i=1,\cdots,n}\kappa_i$ and
$K:=\sum_{i=1}^n \sup_{\theta_i \in \Theta_i}\|\xi_d(\theta_{i})\|$.
%Moreover,
Thus
\bgeqn
\label{eq:uniform-conver-3}
\sup_{V\in \mathcal{U}}|\phi^M(V,V^*)-\phi^M(V,V^M)| \leq \kappa
(n-1)
 K  \|V^*-V^M\|.
\edeqn
Combining (\ref{eq:uniform-conver-2}) and (\ref{eq:uniform-conver-3}), we have
\bgeqn
\label{eq:unform-limit}
\sup_{V\in \mathcal{U}}|\phi(V,V^*)-\phi^M(V,V^M)| \leq L\dd_K(\eta,\eta^M) + \kappa
(n-1)
 K  \|V^*-V^M\|.
\edeqn
%as $M \to 0$.
Since $V^M \in \arg\max_{U \in \mathcal{V}_d^M} \phi^M(U,V^M)$,
then
%for any $U^M \in \mathcal{V}_d^M$, we have
\bgeqn
\label{eq:optimal-V-M}
\phi^M(U^M,V^M) \leq \phi^M(V^M,V^M), \forall U^M \in \mathcal{V}_d^M.
\edeqn
Without loss of generality, we assume that $U^M \to U$.
By Proposition \ref{p:V-cover}, $U \in \mathcal{V}_d$.
Moreover, by (\ref{eq:unform-limit}) and continuity of $\phi$, we have
\bgeqn
|\phi^M(V^M,V^M)-\phi(V^*,V^*)|
&=&|\phi^M(V^M,V^M)-\phi(V^M,V^*)\nonumber\\
&&+\phi(V^M,V^*)-\phi(V^*,V^*)|
%\nonumber\\
%&\leq& 2 (L\dd_K(\eta,\eta^M) + \kappa
%(n-1)
% K  \|V^*-V^M\|).
 \to 0,
 \label{eq:optimal-V-M-1}
\edeqn
%and
Likewise
\bgeqn
|\phi^M(U^M,V^M)-\phi(U,V^*)|
&=&|\phi^M(U^M,V^M)-\phi(U^M,V^*)\\\nonumber
&&+\phi(U^M,V^*)-\phi(U,V^*)| \to 0.
\label{eq:optimal-V-M-2}
\edeqn
Hence by driving $M \to \infty$, we obtain from
(\ref{eq:optimal-V-M})-(\ref{eq:optimal-V-M-2}) that
$$
\phi(U,V^*) \leq \phi(V^*,V^*).
$$
Since $U$ can be arbitrary in $\mathcal{V}_d$ (because Proposition \ref{p:V-cover} ensures that for any $U\in \mathcal{V}_d$, we can find
$U^M\in \mathcal{V}_d^M$ such that $U^M\to U$),
 we arrive at
$$
V^* \in \inmat{argmax}_{U \in \mathcal{V}_d} \phi(U,V^*),
$$
which implies $V^*$ is an equilibrium of (BNE-app$'$-dis).

{\color{black}
%\textbf{Keep it for the time being}
Part (ii).
%For the second part,
By using Lemma \ref{l-stability}, we have
\bgeq
\mathbb{D}(V^M, {\cal V}^*)
&\leq& \left\{\frac{3}{\alpha} \sup_{V \in  {\cal U}} |\phi(V,V^*)-\phi^M(V,V^M)|\right\}^{\frac{1}{\nu}}\\
&\leq& \left\{\frac{3}{\alpha} (L \dd_K(\eta,\eta^M) + \kappa (n-1) K\|V^*-V^M\|)\right\}^{\frac{1}{\nu}} \quad \inmat{(using (\ref{eq:unform-limit}))}.
\edeq
When ${\cal V}^*=\{V^*\}$ and  $M$ is sufficiently large  with
$\|V^*-V^M\|\leq \left(\frac{\alpha}{4\kappa(n-1) K}\right)^{\frac{1}{1-\nu}}$,
 we  obtain
\bgeqn
\|V^*-V^M\| \leq \left( \frac{12L}{\alpha} \dd_K(\eta,\eta^M)\right)^{\frac{1}{\nu}}.
\label{eq:error-bnd-V*-V^M}
\edeqn

Part (iii). In this case, $\phi^M(V,V^*)$ is an ordinary sample average approximation of $\phi(V,V^*)$.
Since $\Theta$ is compact and $u_i(v_i^T \xi_d(\theta_i),(v_{-i}^*)^T\xi_d(\theta_{-i}),\theta)$ is Lipschitz continuous in $v_i$ for $i=1,\cdots,M$, we may use \cite[Theorem 5.1]{SX08}
to establish uniform exponential convergence of
$\phi^M(V,V^*)$ to $\phi(V,V^*)$.
From the proof of Part (ii), we know
{\color{black}
\bgeq
\|V^M- V^*\|^{\nu}
&\leq& \frac{3}{\alpha} \sup_{V \in  {\cal U}} |\phi(V,V^*)-\phi^M(V,V^M)|\\
&\leq& \frac{3}{\alpha} \sup_{V \in  {\cal U}} |\phi(V,V^*)-\phi^M(V,V^*)|+ \frac{3}{\alpha} \sup_{V \in  {\cal U}} |\phi^M(V,V^*)-\phi^M(V,V^M)|\\
&\leq& \frac{3}{\alpha} \sup_{V \in  {\cal U}} |\phi(V,V^*)-\phi^M(V,V^*)|+  \frac{3}{\alpha}\kappa(n-1) K  \|V^*-V^M\|\quad \inmat{(using (\ref{eq:uniform-conver-3}))}\\
&\leq& \frac{3}{\alpha} \sup_{V \in  {\cal U}} |\phi(V,V^*)-\phi^M(V,V^*)|+ \frac{3}{4}  \|V^*-V^M\|^{\nu}\quad \left(\inmat{when $\|V^*-V^M\|\leq \left(\frac{\alpha}{4\kappa(n-1) K}\right)^{\frac{1}{1-\nu}}$}\right).
\edeq
%and furthermore, we obtain
Consequently, we have
\bgeq
\|V^M- V^*\| \leq \left( \frac{12}{\alpha} \sup_{V \in  {\cal U}} |\phi(V,V^*)-\phi^M(V,V^*)|\right)^{\frac{1}{\nu}}.
\edeq
The rest follows from exponential convergence of
$\sup_{V \in  {\cal U}} |\phi(V,V^*)-\phi^M(V,V^*)|$ via \cite[Theorem~5.1]{SX08}.
 We omit the details.
%The proof is complete.
}
}
\hfill $\Box$

{\color{black}
It might be helpful to add some comments on the theorem. Part (i) is a kind of qualitative
convergence statement which guarantees the convergence but is short of giving the rate of convergence.
Part (ii) strengthens the result by giving an explicit error bound for $V^M$ under some growth condition.
%  $C(\dd_K(\eta,\eta^M)^{\frac{1}{\nu}}$.
It is important to note that $\dd_K(\eta,\eta^M)^{\frac{1}{\nu}}$ depends on the dimension of $\theta$ which means when the dimension is high, the bound could be rough, in other words, it is subject to curse of dimensionality.
Part (iii) addresses this issue, it says that when the points of $\tilde{\Theta}^M$ are given through iid samples,
the probability of $V^M$ deviating from $V^*$ reduces at exponential rate with increase of the sample size
when
$V^M$ falls into
$\left(\frac{\alpha}{4\kappa(n-1) K}\right)^{\frac{1}{1-\nu}}$-range of $V^*$ with probability $1$.
From \cite[Theorem 5.1]{SX08}, one can see that the constants $C(\varepsilon)$ and $\lambda(\varepsilon)$
depend on the dimension of $V$ and the size of its domain rather than the dimension of $\theta$.
This means when the dimension of $\theta$ is low, we may use the optimal quantization method
whereas when the dimension of $\theta$ is high, it might be more efficient to use the well-known sample average approximation method.
}

\section{Applications to rent-seeking contests}

%{\color{red}
%1. Review on rent-seeking;
%
%%(1) Tullock, G. (1980). Efficient rent seeking. In J. Buchanan, R. Tollison \& G. Tullock (Eds.), Toward a theory
%%of the rent-seeking society (pp. 97–112). College Station: Texas A \& M Univ. Press.
%
%(2) Ewerhart, C., Quartieri, F., 2015. Unique equilibrium in incomplete information contests with budget constraints. Working paper.
%
%Ewerhart and Quartieri
%considered a more general class of rent-seeking games and
%obtain a sufficient condition for the existence of a unique Bayesian Nash equilibrium,
%which is also based upon strict monotonicity of the payoff gradient.
%
%
%
%
%
%2. extend $g_i(y_i)=y_i$ to general case; two paper extends to three paper symmetric  extends to asymmetric
%}

In this section, we apply the theory on existence of equilibria and approximation schemes established in the previous sections
to rent-seeking contests
with incomplete information.
%A rent-seeking contest is a situation where players
%expend costly effort to gain a reward. Many conflict situations can be described by  rent-seeking contests including political
%campaigns, patent races, war fighting, lobbying efforts, labor market competition, legal battles and professional sports \cite{Fey08}.
There have been extensive literature on Bayesian behavioural function equilibrium or PSNE for studying
% equilibria arising from
such contests.
For instance,
%{\color{red}
%Tullock \cite{Tull80} is the first to investigate this contest.
%}
%{\color{blue} Shaoyan:
%The reference does not seem to be correct, it should be
%Tullock,G.,1967.Thewelfarecostoftariffs,monopoliesandtheft.West.Econ.J.5, 224–232.
%Please double check.
%}
Fey \cite{Fey08} considered rent-seeking contests with two players where each player has private information about his own
cost of effort and
modelled
 them as a Bayesian game where each player's cost is drawn from a distribution of possible costs.
He investigated existence of equilibria for the cases when the distribution of costs is
 discrete and continuous.
Ewerhart \cite{EC14} advanced the research by showing existence of unique PSNE where the contest success function is of logit form with concave impact functions and each player's private information may
relate to either costs or valuations.
%{\color{red}
Ewerhart and Quartieri \cite{EQ15}
considered a more general class of rent-seeking contests and
obtained a sufficient condition for the existence of a unique Bayesian Nash equilibrium.
%}
Here we
follow primarily
%Ewerhart's model.
%{\color{red}
Fey's model.
%}

%Suppose that there are
Consider a rent-seeking contest with $n$ players ($n\geq 2$)
% players
who
%need
aim to choose a level of costly effort $a_i \geq 0$
in order
to obtain a share of a prize,
and each player's cost
is a linear function of his effort
parameterized by
$\theta_i$.
The value of $\theta_i$ is drawn from a probability distribution $\eta_i$ which is
absolutely continuous with respect to the Lebesgue measure
over its
%with
 support set $\Theta_i=[\alpha_i,\beta_i]$.
%Let
Assuming that $\theta_1,\ldots, \theta_n$
are
 independent,
%Under this assumption,
we can write down
the expected utility of player $i$
%'s with cost $\theta_i$,
%when the opponents play strategy $f_{-i}$,
%is given by
%can be written
as
\bgeqn
\label{eq:rent-seeking-model}
\rho_i(a_i,f_{-i},\theta_i):
=-a_i \theta_i +\int_{\Theta_{-i}} \frac{a_i}{a_i+\sum_{j\neq i} f_j (\theta_j)} d \eta_{-i}(\theta_{-i})\,\inmat{for}\,i\in N,
%=1,\ldots,n,
\edeqn
where
%the utility function takes the form:
\bgeq
u_i(a,\theta)=-a_i \theta_i +\frac{a_i}{a_i+\sum_{j\neq i} a_j}.
\edeq
In the case when
$a_i+\sum_{j\neq i} f_j (\theta_j)=0$,
%$a_i=0$ and $f_j (\theta_j)=0$ for
%Let $\rho_i(a_i,f_{-i},\theta_i)=0$ for $a_i=0$ and $f_{-i}(\theta_{-i})=0$.
%Here, t
$\frac{a_i}{a_i+\sum_{j\neq i} f_j (\theta_j)}$ is set $\frac{1}{n}$, which means
each player gets a fair $\frac{1}{N}$ if no one makes a positive effort.

Note that by letting $a_i=0$,
player $i$ can obtain a payoff of zero,
hence the optimal choice of $a_i$ must satisfy $a_i\geq 0$ with
$\rho_i(a_i,f_{-i},\theta_i) \geq 0$ for all $\theta_i \in \Theta_i$.
Moreover, since the integral in (\ref{eq:rent-seeking-model}) is bounded by $1$,
%this implies that
the optimal choice of $a_i$
must satisfy
$$
0 \leq a_i \leq 1/\theta_i.
% \leq 1/ {\alpha_i}.
$$
Since $\alpha_i =\min  \{\theta_i: \theta_i\in \Theta_i\}$,
%{\color{red} should min}
then we can
 define the action space of player $i$ by
 $\mathcal{A}_i:=[0, 1/ {\alpha_i}]$.

 To fit the problem entirely into our framework, we make the action space of each player a bit more restrictive
 by considering $\mathcal{A}_i^{\epsilon}=[\epsilon, 1/ {\alpha_i}]$ for some small positive constant $\epsilon$.
This is justified in the case when $n=2$. To see this, if player $i$ observes his opponent makes zero effort ($f_{-i}=0$),
then he would clearly be better off by making a small effort $\epsilon$ (the smaller the better but not equal to zero) with expected profit close to $1$. On the other hand, if his rival (player $-i$) sees $i$ plays $\epsilon$, he would
be better off by setting $a_{-i}$ to $\epsilon+\delta$ where $\delta$ is a small positive number of scale, i.e., of scale $\epsilon$.
Moreover, each player would be better off by keeping its
%component
opponent  making positive effort.
%This has been observed by
Fey \cite{Fey08} observed this in the symmetric case
and
%he
asserted that the observation applies to symmetric multi-player case.
As far as we are concerned, it is unclear
whether this is correct or not in asymmetric situations.

With $\mathcal{A}_i^{\epsilon}$ being defined as above, we
%first
consider the following (BNE): find an $n$-tuple of behavioural functions $(f_1^\epsilon,\cdots, f_n^\epsilon)$ such that
\bgeqn
\label{eq:rent-BNE}
  \quad f_i^{\epsilon}(\theta_i) \in\arg \max_{a_i\in \mathcal{A}_i^{\epsilon}} \rho_i(a_i,f_{-i}^{\epsilon},\theta_i), \;
  \forall \theta_i\in \Theta_i\,
  %\inmat{and} \;  i\in N,
\edeqn
for  $i\in N$,
or
equivalently
\bgeq
\rho_i(f_i^{\epsilon}(\theta_i),f_{-i}^{\epsilon},\theta_i) \geq \rho_i(g_i(\theta_i),f_{-i}^{\epsilon},\theta_i), \forall  \theta_i \in \Theta_i,
\edeq
for
%every
any $g_i: \Theta_i \to \mathcal{A}_i^{\epsilon}$, where
$\rho_i$ is defined as in (\ref{eq:rent-seeking-model}).
% and $i\in N$.
Existence of Bayesian behavioural function equilibrium is established by Ewerhart \cite[Lemma 3.1]{EC14}.
Here we show that the existence can also be verified by our theoretical results in Section 3.
%In order to apply the existence result in Section 3 to this model,
To see this, it suffices to
%we  first need to
verify the conditions in Assumption \ref{A:BNE-Main}.

%It is easy to obtain that
From the definition of the problem, we observe that: (a) $u_i$ is continuous over $\mathcal{A}^{\epsilon} \times \Theta$;
(b) $\eta_i(S| \theta_i)=\eta_i(S)$ for any measurable set $S \subset \Theta_{-i}$,
which implies $\eta_i(S| \theta_i)$ is continuous in $\theta_i$;
(c) for any $\theta_i, \theta_i'\in \Theta_i$,
$$
|\rho_i(a_i,f_{-i}^{\epsilon},\theta_i)-\rho_i(a_i,f_{-i}^{\epsilon},\theta_i')|
= a_i |\theta_i -\theta_i'| \leq 1/ {\alpha_i}|\theta_i -\theta_i'|.
$$
%Now we only need to prove the strong concavity of $u_i$ over $\mathcal{A}_i^{\epsilon}$.
Moreover,
%by simple calculations,
%we have
\bgeq
(u_i)'_{a_i}(a,\theta)=-\theta_i+\frac{\sum_{j\neq i} a_j}{(a_i+\sum_{j\neq i} a_j)^2}
\edeq
and
\bgeq
(u_i)''_{a_i}(a,\theta)=\frac{-2\sum_{j\neq i} a_j}{(a_i+\sum_{j\neq i} a_j)^3}.
\edeq
Since $a_i \in \mathcal{A}_i^{\epsilon}=[\epsilon,1/\alpha_i]$,
\bgeqn
\label{eq:rent-concave}
(u_i)''_{a_i}(a,\theta) \leq \frac{-2(n-1) \epsilon}{ (\sum_{i=1}^n 1/\alpha_i)^3}<0,
\edeqn
%and hence
which means that $u_i$ is strongly concave over $\mathcal{A}_i^{\epsilon}$.
Therefore, all the conditions of Assumption \ref{A:BNE-Main}
are satisfied here.
%fulfilled.
By Theorem \ref{t-NE-existence},
the problem (\ref{eq:rent-BNE}) has a continuous behavioural function equilibrium
denoted by $(f_1^{\epsilon},\ldots,f_n^{\epsilon})$.

%Now, w
We now move on to
%e will verify
%problem (\ref{eq:rent-BNE}) has an unique equilibrium.
show uniqueness of the equilibrium.
%Moreover, i
It is easy to observe that
%From (\ref{eq:rent-concave}),
%we know that
%$u_i(a,\theta)$ is strictly concave in $a_i$.
%Since
%$$
%(u_i)'_{a_{j}}(a,\theta)=\frac{-1}{(a_i+\sum_{j\neq i} a_j)^2}\,\inmat{and}\,
%(u_i)''_{a_{j}}(a,\theta)=\frac{2}{(a_i+\sum_{j\neq i} a_j)^3}>0
%$$
%for $j\neq i$, then
$u_i$ is convex in $a_{-i}$ in that
the second term of $u_i$ can be viewed as
composition of a convex function and a linear function of $a_{-i}$ .
%Furthermore,
Moreover,
%due to
since $\sum_{i=1}^n u_i(a,\theta)= 1-\sum_{i=1}^n a_i \theta_i$,
$\sum_{i=1}^n u_i(a,\theta)$ is concave in $a$.
%Hence it follows from
By \cite[Lemma 5]{Ui16},
% that
condition (\ref{eq:H-monotone}) is satisfied here and hence
%,and then
the uniqueness
%of equilibrium is guaranteed by
follows by Theorem \ref{t-unique-BNE}.

Let $\{\epsilon_M\}$ be a sequence of  small positive
numbers which
monotonically decreases
%ing
% numpositive $\epsilon_M \searrow 0$ as
to zero
as $M \to \infty$.
We consider convergence of the corresponding
 behavioural function equilibrium
$f^M:=(f_1^{\epsilon_M},\ldots,f_n^{\epsilon_M})$.
% be a BNE  of (\ref{eq:rent-BNE}) for $\epsilon=\epsilon_M$.
In
%Following
\cite[Lemma 3.2 and Theorem 3.4]{EC14}, Ewerhart showed that
$f^M$ has a uniformly converging subsequence which
 converges  to  a continuous behavioural function $f^*$ and
 $f^*$ is a PSNE. Here we draw a slightly stronger conclusion
 by showing $f^*$ is indeed a continuous behavioural function equilibrium
 in the unconstrained contest (an equilibrium of (\ref{eq:rent-BNE}) with $\epsilon=0$).

%$f^*$ has the following property
%and $f^*$ is a BNE of (\ref{eq:rent-BNE}) for $\epsilon=0$.
%Similar to  \cite[Lemma 3.3, Theorem 3.4]{EC14},
%We need the following technical result which is strengthened from
%\cite[Lemma 3.3]{EC14}.

\begin{prop}
\label{p:rent-geq0}
Suppose that $\eta_i$ is absolutely continuous w.r.t. the Lebesgue measure over $\Theta_i$ for
 $i\in N$. Then

 \begin{itemize}

 \item[(i)] there exists some player $i$ such that $f_i^*(\theta_i)>0$ for all
$\theta_i \in \Theta_i$;

\item[(ii)]
$f^*$ is a behavioural function equilibrium
%BNE
of problem (\ref{eq:rent-BNE})
% for
with $\epsilon=0$.
\end{itemize}

\end{prop}

\noindent
\textbf{Proof.}
The proof follows
essentially from
%fundamentally
%similar to
similar proofs of \cite[Lemma 3.3 and Theorem 3.4]{EC14}. We include it for completeness.
%the minor difference.
%We prove it by
 For the simplicity of notation,
we assume without loss of  generality
that
%, we assume
$f^M$ converges to $f^*$.

Part (i). Assume for the sake of a contradiction
%Suppose that
that for
%each
any $i\in N$, there exists a point $\theta_i^* \in \Theta_i$ such that
$f_i^*(\theta_i^*)=0$.
By the continuity of $f_i^*$ and uniform convergence of $f^M$ to $f^*$,
for any $\delta>0$, there exist a positive number $M_0$ and
a neighborhood $\B(\theta_i^*)$ of $\theta_i^*$
such  that
$$
f_i^*(\theta_i) \leq \delta\, \quad \inmat{and}\, \quad f_i^M(\theta_i) \leq 2\delta, \quad \forall \theta_i \in \B(\theta_i^*)
$$
when $M \geq M_0$.
On the other hand,
by the first order optimality condition,
%for
%at type $\theta_i^*$,
% in (BNE),
%we have
\bgeq
0 \geq \int_{\B(\theta_{-i}^*)} \frac{\sum_{j\neq i} f_j^M (\theta_j)}{(f_i^M(\theta_i)+\sum_{j\neq i} f_j^M (\theta_j))^2} d \eta_{-i}(\theta_{-i})-\theta_i, \;\; \forall \theta_i\in \B(\theta_i^*).
\edeq
%By i
Integrating w.r.t. $\theta_i$ over $\B(\theta_i^*)$ and summing over $i=1,\ldots,n$,
we obtain
\bgeq
0 \geq \int_{\B} \frac{n-1}{\sum_{i=1}^n f_i^M (\theta_i)} d \eta(\theta)- \sum_{i=1}^n \int_{\B(\theta_i^*)}\theta_i d \eta_i(\theta_i),
\edeq
where $\B=\B(\theta_1^*)\times \cdots \times \B(\theta_n^*)$.
%Due to
Since  $f_i^M(\theta_i)\leq 2\delta$ over $\B$ and $\delta$ can be
%sufficiently
arbitrarily small,
the above inequality does not hold.
%The proof is complete.
%\hfill $\Box$
%
%}
%
%
%\begin{prop}
%$f^*$ is a BNE of (\ref{eq:rent-BNE}) for $\epsilon=0$.
%\end{prop}
%
%\noindent
%\textbf{Proof.}

Part (ii). Since $f^M$ is a behavioural function equilibrium of  (\ref{eq:rent-BNE}) with $\epsilon=\epsilon_M$,
%that is,
\bgeq
\rho_i(f_i^M(\theta_i),f_{-i}^M,\theta_i) \geq \rho_i(a_i,f_{-i}^M,\theta_i), \forall a_i \in \mathcal{A}_i^{\epsilon_M},
\edeq
for each fixed $\theta_i \in \Theta_i$.
%a_i \in \mathcal{A}_i^{\epsilon_M}$ and $i\in N$.
Therefore, if $f_{-i}^*>0$ with positive measure, then we may drive
%taking
$M$ to infinity and
% \to \infty$, we have
obtain
\bgeqn
\label{eq:rent-0}
\rho_i(f_i^*(\theta_i),f_{-i}^*,\theta_i) \geq \rho_i(a_i,f_{-i}^*,\theta_i), \forall  \theta_i \in \Theta_i
\edeqn
for $a_i \in \mathcal{A}_i$ and $i\in N$.
If $f_{-i}^*(\theta_{-i})=0$ for all $\theta_{-i}\in \Theta_{-i}$,
then by
%Lemma \ref{l:rent-geq0}
Part (i),
%that
$f_i^*(\theta_i)>0$ for all $\theta_i \in \Theta_i$.
As we commented before, the optimal strategy for the player is to set $a_i$
close to $0$ but not $zero$ in which case the expected profit would be close to $1$ as opposed to $1/N$
with $a_i=0$.
%The expect profit
%Then we also have
%In that case,
This shows (\ref{eq:rent-0}) still  holds for all $a_i\geq 0$.
%Note that
%%Since
% $\rho_i(a_i,f_{-i}^*,\theta_i)$ is lower semicontinuous in $a_i$,
% so
%(\ref{eq:rent-0})
%%also
%holds
%%for
%at $a_i  = 0$.
%%{\color{red}
%%Shaoyan, I cannot see why we need this as $\rho_i(0,f_{-i}^*,\theta_i)=0$.}
%% for this case.
The proof is complete.
\hfill $\Box$

The weakness of this proposition is that we are short of claiming whether or not every cluster point of
$\{f^M\}$
is a continuous behavioural function equilibrium
 of the unconstrained contest when the sequence has multiple cluster points.

For the case that $\epsilon=0$, we can obtain from \cite[Theorem 1]{Fey08} that
(\ref{eq:rent-BNE}) has an equilibrium with continuous behavioural functions.
Fey \cite{Fey08} proposed a standard iterative method
%to
for solving (\ref{eq:rent-BNE}).
Here we apply the polynomial decision rule and discretization scheme  discussed in Section 4 to
model (\ref{eq:rent-BNE})
%for
with $\epsilon=0$,
which means we solve the following (BNE):
\bgeqn
\label{eq:discrete-rent}
V^M \in \arg\max_{V \in \mathcal{V}_d^M} \sum_{i=1}^n
\bbe_{\eta^M}\left[
-v_i^T \xi_d(\theta_i) \theta_i +
\frac{v_i^T \xi_d(\theta_i)}{v_i^T \xi_d(\theta_i)+\sum_{j\neq i} (v_{j}^M)^T \xi_d(\theta_{j})}
\right],
\edeqn
where
\bgeq
(\mathcal{V}_d^M)^i :=\left\{v_i \in \mathbb{R}^{d+1}: 0 \leq v_i^T \xi_d(\theta_i) \leq 1/{\alpha_i}, \, \forall \theta_i \in \tilde{\Theta}_i^M\right\}.
\edeq

{\color{black}
We have carried out numerical tests
on a
%couple of
symmetric rent-seeking contest and an asymmetric
one.
%Note that f
In the symmetric case, Fey \cite{Fey08} obtained a
numerical solution by discretizing
the problem over a grid of $100$ elements and believed
that it provides a good approximation to the unknown true equilibrium.
  Moreover,  his focus was
  %to see
 on the difference between
 the equilibrium effort level under incomplete information
% differs from
and the effort level under complete information.
 Ewerhart \cite{EC14} considered the asymmetric case.
Here our focus is how the solution obtained through our approximation schemes
is affected by variation of order of the polynomials $d$
and sample size $M$.

%In the following examples, we consider
In the numerical experiments, we
concentrate on two player games
%the number of players is 2
and use a heuristic method, the Gauss-Seidel-type Method \cite{FHS14},
to solve problem (\ref{eq:discrete-rent}).
The tests are carried out in MATLAB 7.10.0 installed on a Dell-PC with
Windows 10 operating system and Intel Core i3-2120 processor.

 %see \cite{FHS14} for more detail.

\begin{alg}\rm
\label{alg1}
Let $V^0=(v_1^0, v_2^0)$ and set $k=0$.

\begin{itemize}

\item[\inmat{Step 1.}] For given $V^k=(v_1^k,v_2^k)$, solve
\bgeqn
\begin{array}{cl}
\displaystyle \max_{v_2 \in \R^{d+1}}   & \sum_{j=1}^M {p_j} \left(
-v_2^T \xi_d(\theta_2^j) \theta_2^j +
\frac{v_2^T \xi_d(\theta_2^j)}{v_2^T \xi_d(\theta_2^j)+ (v_1^k)^T \xi_d(\theta_1^j)}
\right)\\
\inmat{s.t.}               &0 \leq v_2^T \xi_d(\theta_2^j) \leq 100,\, j=1,\ldots,M.
\end{array}
\edeqn
Let $v_2^{k+1}$ denote the optimal solution.
Then solve
\bgeqn
\begin{array}{cl}
\displaystyle \max_{v_1 \in \R^{d+1}}   & \sum_{j=1}^M {p_j} \left(
-v_1^T \xi_d(\theta_1^j) \theta_1^j +
\frac{v_1^T \xi_d(\theta_1^j)}{v_1^T \xi_d(\theta_1^j)+ (v_2^{k+1})^T \xi_d(\theta_2^j)}
\right)\\
\inmat{s.t.}               &0 \leq v_1^T \xi_d(\theta_1^j) \leq 100,\, j=1,\ldots,M.
\end{array}
\edeqn
Let $v_1^{k+1}$ denote the optimal solution.

\item[\inmat{Step 2.}]
If $V^{k+1} =V^k$,  stop. Otherwise,
let $V^k :={V}^{k+1}$, go to Step 1.

\end{itemize}
\end{alg}

\begin{ex}[Symmetric rent-seeking contests \cite{Fey08}]
\rm
Let $n=2$, $\Theta_1=\Theta_2=[0.01, 1.01]$,
and $\mathcal{A}_1=\mathcal{A}_2=[0, 100]$.
Suppose that $\theta_1$ and $\theta_2$ are independent
and uniformly distributed over $\Theta_1$ and $\Theta_2$ respectively.

In order to look into the performance of our approximation schemes,
we have carried out two sets of experiments
%under variation of
with respect to change  of the order of the polynomials $d$
and the sample size $M$.
%One is f
We start with
%For
fixed sample size $M$ and
%we
%test
investigate the
performance of
the approximate
%schemes
behavioural function equilibrium as
%order of the polynomials
 $d$
increases.
Figure \ref{Fig1} visualizes
changes of the behavioural functions of both players (they are identical as the game is symmetric),
 we can see that for fixed sample size $M=4900$,
% {\color{red} should it be 4900?}
there are
sizable
%a big
%change
shifts of the
behavioural function curves over the interval $[0.01,0.2]$
%in the graphs
as
%order of the polynomials
$d$ increases from $5$ to $8$
%and
%the graphs is nearly the same
but
%no
%significant difference
stabilizes after $d$ reaches $8$.
%for
%between $d=8$ and $d=9$.
% over the interval $[0.01,1.01]$.

We then move on to
%The other is for the fixed order of the polynomials $d$,
%we
examine
%the
%performance with different
the impact on the approximate behavioural functions at equilibrium
%against variation of
as the sample size $M$ increases for fixed $d$.
%In this example, since the set $\Theta_1$ and $\Theta_2$ are identical,
%the behavioural functions of player 1 and player 2 are the same.
%The results are depicted in Figures \ref{Fig1} and \ref{Fig2} for player 1.
Figure \ref{Fig2} displays
%the
%performance schemes
changes of the behavioural functions at equilibrium
%where we fix up the order of the polynomials $d=9$
%but
%and
when we change the sample size from $100$ to $400, 900$ and $1600$ with fixed $d=9$.
It can be seen from Figure \ref{Fig2} that after $M$ reaches
$900$, there is no significant change.
% and $M=1600$,
%the corresponding graphs are similar.

\begin{figure}[!htp]
 \centering                \includegraphics[width=17cm,height=8cm]{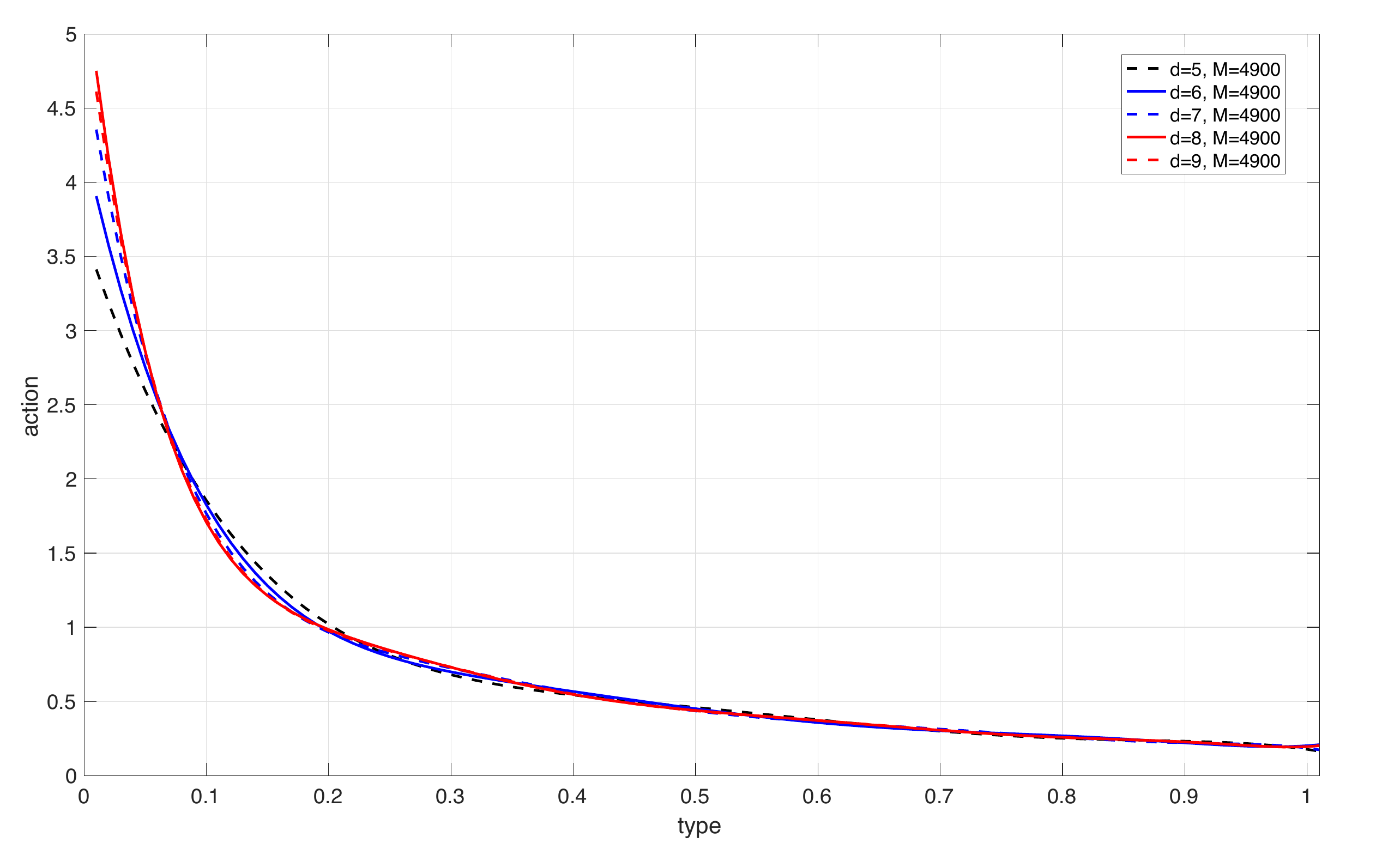}
      \caption{Performance v.s.  order of the polynomial $d$,  Example \ref{ex1}.}
              \label{Fig1}
 \end{figure}

\begin{figure}[!htp]
 \centering                \includegraphics[width=17cm,height=8cm]{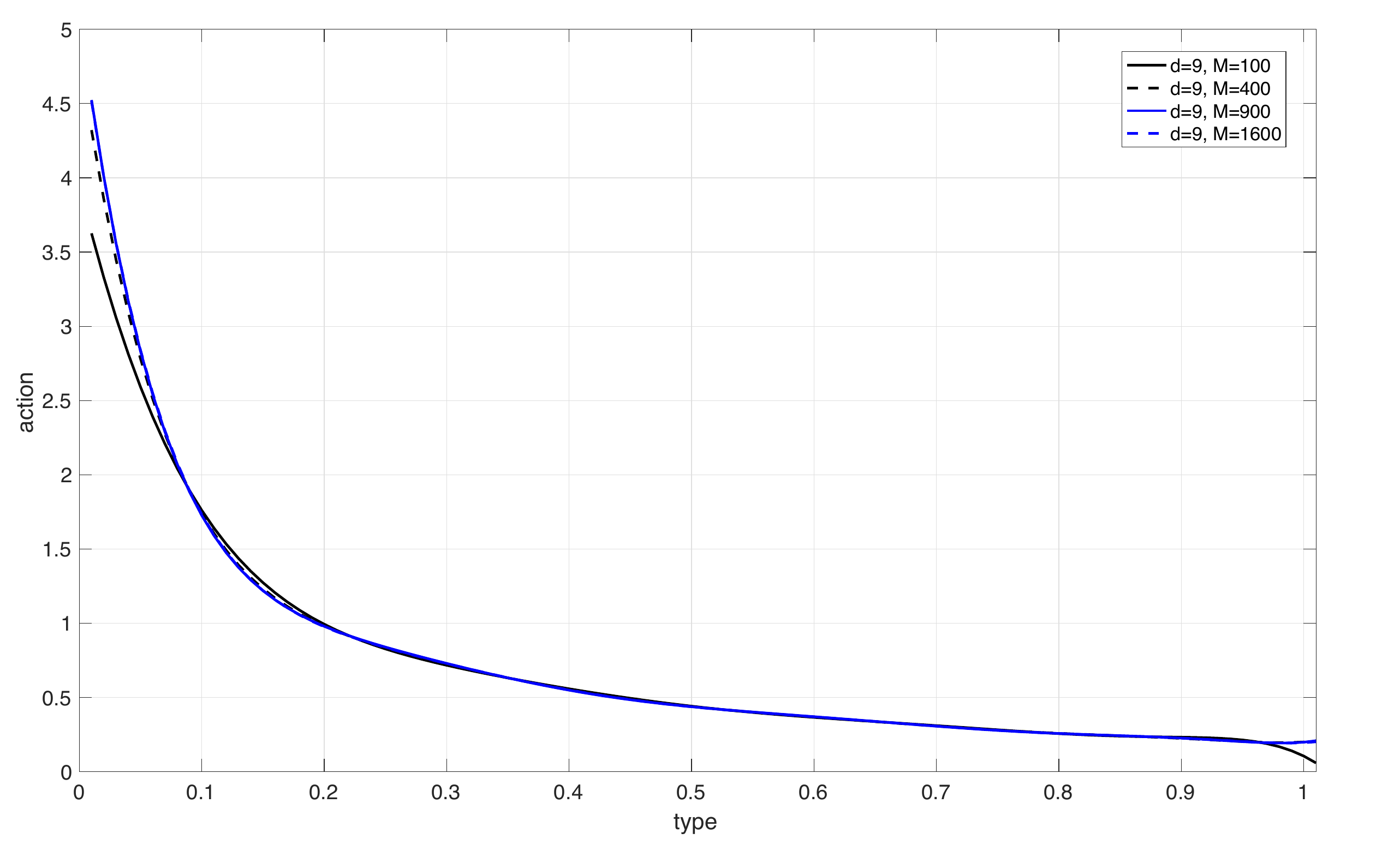}
      \caption{Performance v.s.  sample size $M$,  Example \ref{ex1}.}
              \label{Fig2}
 \end{figure}
\label{ex1}

Throughout the experiments, the samples are chosen by the discretization scheme
discussed in Section 4. For example, in the first set of experiments,
we pick up $70$ points evenly spread over $\Theta_1$ and $\Theta_2$
respectively and use them to form $4900$ grid points over the space of $\Theta_1\times \Theta_2$.
So these samples are generated in a deterministic manner.
\label{ex-sym-rent-seek}
\end{ex}

Next, we examine the approximation scheme by applying it to an asymmetric rent-seeking contest
with both players having identical expected utility functions and action spaces as in Example \ref{ex-sym-rent-seek} but
with different type sets.

%In the following example, we consider the case that $\Theta_1$ and $\Theta_2$ are different.

\begin{ex}[Asymmetric rent-seeking contests]\rm
%{\color{red} (please give a reference. The type sets are different with \cite{EC14})}
Let $n=2$,
$\Theta_1=[0.01, 1.01]$,
$\Theta_2=[0.01, 2.01]$
and $\mathcal{A}_1=\mathcal{A}_2=[0, 100]$.
As in Example \ref{ex-sym-rent-seek}, we assume that $\theta_1$ and $\theta_2$ are independent
and uniformly distributed over $\Theta_1$ and $\Theta_2$ respectively.
This example is varied from  Ewerhart  \cite{EC14} where $\Theta_2=[0.51, 5.51]$ whereas
all other settings are the same.

We
have carried out two sets of experiments as
in Example \ref{ex1}.
The results are depicted in Figures \ref{Fig3} and \ref{Fig4}.
Figure \ref{Fig3}
visualizes changes of the approximate behavioural functions at equilibrium
 for player 1 and player 2
when the order of the polynomials $d$
increases from $5$ to $8$ with fixed sample size $M=5000$.
Figure \ref{Fig4}
depicts changes of the approximate behavioural functions at equilibrium
when the sample size $M$ changes from $72$ to $392$ with fixed order of the polynomials $d=8$.
{\color{black}
Note that different from Example \ref{ex1},
the size of interval $\Theta_2$ is twice of $\Theta_1$,
so we pick up $K$ points and $2K$ points evenly from
 $\Theta_1 $ and $\Theta_2$ with $K=6, 10, 13, 14$
 and use them to generate  $K \times 2K$ grid points/samples.
% when we use discretization scheme in Section 4, e.g.
%$72=6 \times12$ and $392=14 \times 28$.
}

%{\color{red} (Shaoyan, please explain how the number 72 and 392 are obtained.)}

%{\color{blue}
%The approximation schemes  
The preliminary numerical tests show that our approximation schemes work very well.
Note that it is possible to reformulate problem (\ref{eq:discrete-rent})
into a nonlinear complementarity problem (NCP) through first order optimality conditions and consequently we may
replace Algorithm \ref{alg1} with an existing NCP solver such as PATH.
Since the reformulation is equivalent, it does not affect the test results but may avoid the iterative
process.
%}

\begin{figure}[!h]
 \centering                \includegraphics[width=17cm,height=8cm]{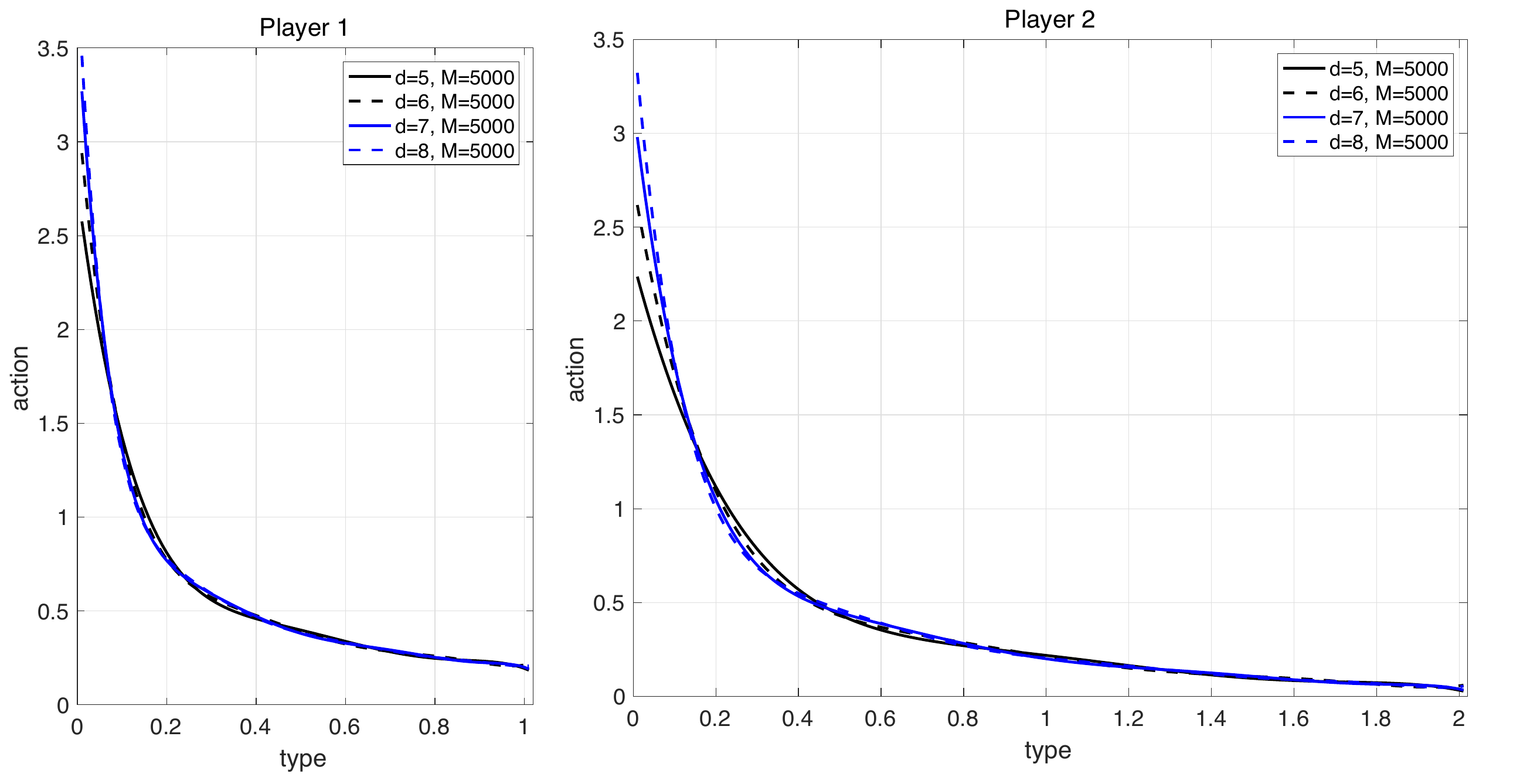}
      \caption{Performance v.s.  order of the polynomial $d$,  Example \ref{ex2}.}
              \label{Fig3}
 \end{figure}

\begin{figure}[!h]
 \centering                \includegraphics[width=17cm,height=8cm]{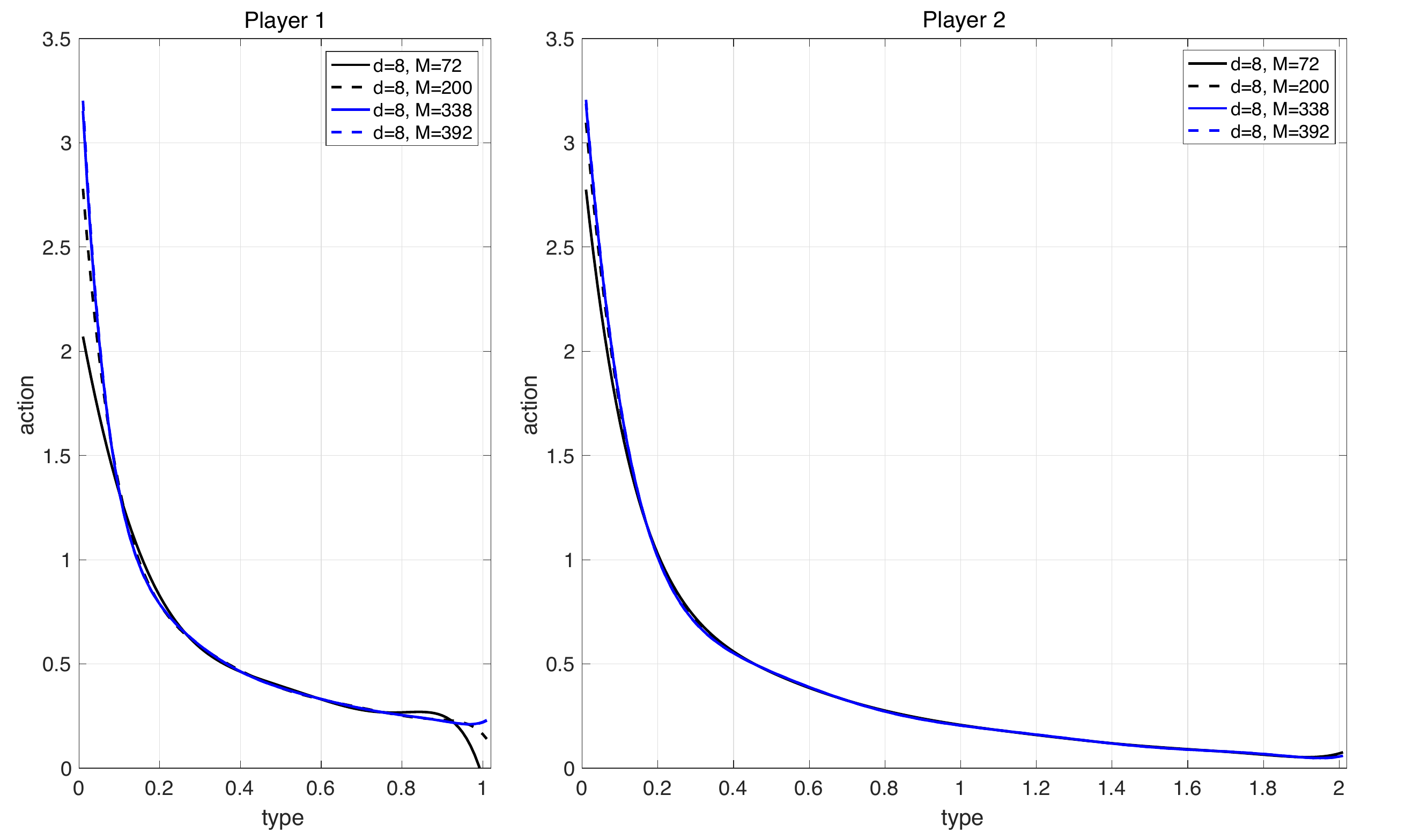}
      \caption{Performance v.s.  sample size $M$,  Example \ref{ex2}.}
              \label{Fig4}
 \end{figure}
\label{ex2}
\end{ex}
}

\end{document}